 \documentclass[final,5p,times,twocolumn]{elsarticle}

\usepackage{amsmath,amsfonts,amssymb,mathtools}
\usepackage{bbm}
\usepackage{graphicx}
\usepackage{exscale}

\biboptions{sort&compress}

\usepackage{clrscode}
\usepackage{array}
\usepackage{hhline}
\usepackage{color, colortbl}
\definecolor{myGray}{rgb}{0.95,0.95,0.95}

\newcommand{\E}[1]{\mathop{{\rm \bf E}\!\left\{#1\right\}}\nolimits}
\newcommand{\var}[1]{\mathop{{\rm \bf var}\!\left\{#1\right\}}\nolimits}

\newtheorem{lmm}{Lemma}
\newproof{pf}{Proof}

\newdefinition{rmk}{Remark}
\newdefinition{dfn}{Definition}
\newtheorem{exmp}{Example}

\begin{document}

\begin{frontmatter}
\title{Continuous-discrete derivative-free extended Kalman filter based on Euler-Maruyama and It\^{o}-Taylor discretizations: Conventional and square-root implementations}

\author[CEMAT]{Maria V. Kulikova\corref{cor}} \ead{maria.kulikova@ist.utl.pt} \cortext[cor]{Corresponding
author.}

\author[CEMAT]{Gennady Yu. Kulikov} \ead{gennady.kulikov@tecnico.ulisboa.pt}

\address[CEMAT]{CEMAT, Instituto Superior T\'ecnico, Universidade de Lisboa, Av.~Rovisco Pais, 1049-001 Lisboa, Portugal.}

\begin{abstract}
In this paper, we continue to study the derivative-free extended Kalman filtering (DF-EKF) framework for state estimation of continuous-discrete nonlinear stochastic systems. Having considered the Euler-Maruyama and It\^{o}-Taylor discretization schemes for solving stochastic differential equations, we derive the related  filters' moment equations based on the derivative-free EKF principal. In contrast to the recently derived MATLAB-based {\it continuous-discrete} DF-EKF techniques, the novel DF-EKF methods preserve an information about the underlying stochastic process and provide the estimation procedure for a fixed number of iterates at the propagation steps. Additionally, the DF-EKF approach is particularly effective for working with stochastic systems with highly nonlinear and/or nondifferentiable drift and observation functions, but the price to be paid is its degraded numerical stability (to roundoff) compared to the standard EKF framework. To eliminate the mentioned pitfall of the derivative-free EKF methodology, we develop the conventional algorithms together with their stable square-root implementation methods. In contrast to the published DF-EKF results, the new square-root techniques are derived within both the Cholesky and singular value decompositions. A performance of the novel filters is demonstrated on a number of numerical tests including well- and ill-conditioned scenarios.
\end{abstract}


\end{frontmatter}

\section{Introduction}\label{sect1}

In recent years, there has been an increasing interest in the development and practical applications of derivative-free Bayesian filtering methods with the Kalman filtering structure, e.g., the Unscented Kalman filtering (UKF) methods proposed in~\cite{Julier1995,Julier2000,Julier2004}, the Cubature Kalman filtering (CKF) algorithms in~\cite{Haykin2009,Haykin2010,JiXi13,Haykin2018}, the Gauss-Hermit Quadrature Kalman filtering (GHQF) techniques in~\cite{Ito2000,Haykin2007} and many others. The derivative-free principal has been proposed for the Extended Kalman filtering (EKF) framework as well. More precisely, the derivative-free EKF (DF-EKF) has been developed for the {\it discrete-time} stochastic systems in~\cite{quine2006derivative}. Recently, the DF-EKF framework has been extended on the {\it continuous-discrete} filtering in~\cite{kulikova2023derivative} where the new methods are based on the MATLAB integrators for solving ordinary differential equations (ODEs). As a result, the suggested {\it continuous-discrete} DF-EKF techniques provide an accurate estimation procedure, but usually at the cost of higher computational demands due to discretization error control involved into adaptive MATLAB ODEs solvers.

An alternative methodology for designing the {\it continuous-discrete} filters' implementation methods is grounded in application of numerical integration methods for solving {\it stochastic differential equations} (SDEs). They are used to discretize the stochastic system at hand and, next, the filters' moment equations are derived for the discretized model. It is important to note that both mentioned strategies for designing the nonlinear Bayesian filters have their own benefits and drawbacks and, in fact, both theoretical approaches are routinely adopted in engineering literature~\cite{Ja70,Frog2012,KuKu14IEEE_TAC,KuKu22Automatica}.

More precisely, the ODE's solvers with an automatic discretization error control utilized at the prediction filtering steps make the filtering methods accurate and simple for implementation because of the built-in fashion of the MATLAB ODEs integrators utilized. They imply an adaptive mesh, which is generated by the solver without any users' work, to satisfy the discretization error bounds and keep the error as small as it is required by users. This makes such methods very useful for solving practical problems with irregular sampling intervals (e.g. when some measurements are missing) and for estimating the so-called stiff systems such as, for instance, the stochastic Van der Pol oscillator example. However, the price to be paid for the mentioned benefits is the increased and unpredictable computational cost. Because of the adaptive mesh generated by the integration schemes according to their discretization error control involved, we do not know in advance how many integration steps will be performed at each sampling interval. This might yield impracticable calculations due to inappropriately large computational time.

In contrast, the filters designed with the use of SDEs numerical integration methods allow to handle the time durations of the filtering methods. Indeed, the SDEs solvers do not allow the discretization error control, but they provide the estimation process for a finite number of steps, which is known and can be set up before the filtering. Thus, such estimation methods are of special interest in applications where the computational time has an essential matter. Additionally, they provide a better approximation of the stochastic terms because of the SDEs integration schemes.  This also means that the information about the stochastic variables at hand is not entirely lost after the numerical integration. In contrast, the ODEs-based filtering approach confines the calculations up to the first two moments, that is, the mean and covariance. If required, any other information about the given stochastic process  is lost and cannot be recovered after the approximation step.

To summarize, the goal of this paper is to develop the {\it continuous-discrete} DF-EKF methods within the Euler-Maruyama and It\^{o}-Taylor discretization schemes applied to the nonlinear stochastic systems. In contrast to the previously published MATLAB-based DF-EKF methods in~\cite{kulikova2023derivative}, the new algorithms keep the information about the underlying stochastic processes and provide the estimation for a fixed number of steps, which is known in advance and can be set up by users. It is anticipated that the It\^{o}-Taylor expansion of strong order 1.5 yields more accurate filtering methods compared to the filters designed under the Euler-Maruyama integrator, which is of strong order 0.5, e.g., see~\cite{KlPl99}. Nevertheless, the simplicity and good estimation quality (when the number of subdivision steps is high enough) are still the attractive features of the Euler-Maruyama-based filtering algorithms and, hence, they deserve some merit. Finally, it should be stressed that the numerical instability to roundoff is a crucial issue for all derivative-free Bayesian filters as discussed in~\cite{Haykin2009,kulikova2023derivative}. To eliminate the mentioned pitfall with respect to the derivative-free EKF methodology, we additionally develop their stable square-root implementation methods. In contrast to the published MATLAB-based DF-EKF algorithms, the new square-root techniques are derived within both the Cholesky factorization and singular value decomposition (SVD). The methods are tested on both the well- and ill-conditioned problems.

The rest of the paper is organized as follows. Section~\ref{sect2} introduces the problem statement and discusses the derivative-free EKF principle for the filter's state and error covariance approximation. Section~\ref{sec:main:conventional} derives the continuous-discrete derivative-free EKF formulas based on the Euler-Maruyama and It\^{o}-Taylor discretization schemes. In Section~\ref{sec:main:squareroot}, we derive numerically stable square-root implementation methods for the novel filters proposed. The results of numerical tests are summarized and discussed in Section~\ref{numerical:experiments}. Section~\ref{Section:conclusion} concludes the paper by the key findings of our research.

\section{Problem statement and derivative-free EKF-based principle for a mean and covariance approximation}\label{sect2}

Consider continuous-discrete stochastic system given by
\begin{align}
dx(t) & = f\bigl(t,x(t)\bigr)dt+Gd\beta(t), \quad t>0,  \label{eq1.1} \\
z_k   & =  h(k,x(t_{k}))+v_k, \quad k =1,2,\ldots \label{eq1.2}
\end{align}
where  $x(t)$ is the $n$-dimensional unknown state vector to be estimated and  $f:\mathbb R\times\mathbb
R^{n}\to\mathbb R^{n} $ is the time-variant drift function and $h:\mathbb R\times\mathbb
R^{n}\to\mathbb R^{m}$ is the observation function. The process uncertainty is represented by the additive noise term where
$G \in \mathbb R^{n\times q}$ is the time-invariant diffusion matrix and $\beta(t)$ is the $q$-dimensional Brownian motion whose increment $d\beta(t)$ is Gaussian white process independent of $x(t)$ and has the covariance $Q\,dt>0$. The $m$-dimensional measurement $z_k = z(t_{k})$ obeys equation~\eqref{eq1.2} and comes at some discrete-time points $t_k$ with the sampling rate (sampling period) $\Delta_k=t_{k}-t_{k-1}$. The measurement noise term $v_k$ in equation~\eqref{eq1.2} is assumed to be a white Gaussian noise with the zero mean and known covariance $R_k>0$, $R_k \in \mathbb R^{m\times m}$. Finally, the initial state $x(t_0)$ and the noise processes are assumed to be statistically independent, and $x(t_0) \sim {\mathcal N}(\bar x_0,\Pi_0)$, $\Pi_0 > 0$.

Let us consider the filtered state estimate and the related error covariance matrix at time instance $t_k$, i.e.
\begin{align*}
\hat x_{k|k} & := \hat x(t_k|t_{k}) = \E{x(t_k)|z_1 \ldots z_{k}} \\
P_{k|k}      & := P(t_k|t_{k}) = \E{[x(t_k) - \hat x_{k|k}][ x(t_k) - \hat x_{k|k} ]^{\top}|z_1 \ldots z_{k}}
\end{align*}
where $\{z_1 \ldots z_{k}\}$ is the data history available for the filter from the measurement device, which is related to the hidden state vector through the measurement model in equation~\eqref{eq1.2}.

Following the derivative-free EKF principle suggested in~\cite{quine2006derivative}, the probabilistic description of the state estimate is approximated by a minimally required number of deterministic vectors. More precisely, $n$ linearly independent vectors are demanded to span the $n$-dimensional vector space. Thus, we consider a set of $n$ discrete vectors $\xi_i$ that approximate the second moment, i.e. the filter error covariance $P_{k|k}$, as follows:
\begin{align}
P_{k|k} & = \frac{1}{n} \sum \limits_{i=1}^{n} \xi_i \xi_i^{\top}. \label{eq:vectors}
\end{align}

Let us define the matrix $[\xi_1 | \ldots | \xi_n]$ collected from vectors $\xi_i$ in such a way that they are located by columns in the discussed matrix. Thus, we may express formula~\eqref{eq:vectors} as follows:
\begin{align}
n P_{k|k} & = \Bigl[\xi_1 \;|\; \ldots \;|\; \xi_n\Bigr] \; \Bigl[\xi_1 \;|\; \ldots \;|\; \xi_n\Bigr]^{\top}, \nonumber \\
\mbox{ i.~e. }   \sqrt{n} P_{k|k}^{1/2} & = \Bigl[\xi_1 \;|\; \ldots \;|\; \xi_n\Bigr]. \label{eq:vectors:1}
\end{align}
where $P_{k|k} = P_{k|k}^{1/2}P_{k|k}^{\top/2}$ with $P_{k|k}^{1/2}$ being the matrix square-root.

The derivative-free EKF principle for the mean and covariance approximation can be formulated as follows~\cite{quine2006derivative}: these $n$-vectors, $\xi_i$, $i=1,\ldots,n$, are utilized to define a set of vectors about the state estimate with the distribution scaled by a constant $\alpha$, i.e. we have
\begin{align}
{\mathcal X}_{i, k|k} & =  \hat x_{k|k} + \frac{1}{\alpha} \xi_i,  & i & = 1, \ldots, n \label{Xi_vectors}
\end{align}
where the scalar parameter $\alpha$ is used to scale the distribution about the mean. It is proved in~\cite{quine2006derivative} that in the limiting case as $\alpha \to \infty$, the derivative-free discrete-time EKF converges to the standard EKF. The choice of the parameter $\alpha$ is adaptable and can be tuned on the application model at hand. More details can be found in the cited paper.

Equation~\eqref{Xi_vectors} can be written in a simple vector-matrix form by introducing the matrix collected from these vectors ${\mathcal X}_{i, k|k} $, $i=1, \ldots,n$, which are located by columns, i.e. we get
\begin{align}
{\mathbb  X}_{k|k}& =\hat x_{k|k} {\mathbf 1}^{\top} + \frac{\sqrt{n}}{\alpha}P_{k|k}^{1/2}, & {\mathbb  X}_{k|k} & =\Bigl[{\mathcal X}_{1,k|k} \;|\; \ldots \;|\; {\mathcal X}_{n,k|k}\Bigr] \label{Xi_vectors:matrix}
\end{align}
where ${\mathbf 1}$ stands for a column vector, whose elements are all equal to one.

Next, taking into account formula~\eqref{Xi_vectors:matrix}, the error covariance matrix, $P_{k|k}$, might be approximated by the vectors ${\mathcal X}_{i, k|k}$, $i = 1, \ldots, n$, in the following way:
{\small
\begin{align}
P_{k|k} & = \frac{1}{n} \sum \limits_{i=1}^{n} \xi_i \xi_i^{\top} = \frac{\alpha^2}{n} \sum \limits_{i=1}^{n} [{\mathcal X}_{i, k|k} - \hat x_{k|k}][{\mathcal X}_{i, k|k} - \hat x_{k|k}]^{\top} = \overline{\mathbb  X}_{k|k} \overline{\mathbb  X}^{\top}_{k|k} \label{P:approximation}
\end{align}
}
where the mean-adjusted and scaled matrix $\overline{\mathbb  X}_{k|k}$ is obtained from the matrix ${\mathbb  X}_{k|k}$ in~\eqref{Xi_vectors:matrix} as follows:

\begin{align}
\overline{\mathbb  X}_{k|k} & = \frac{\alpha}{\sqrt{n}} \bigl[{\mathcal X}_{1,k|k} - \hat x_{k|k}\; |\; \ldots \; | \; {\mathcal X}_{n,k|k} - \hat x_{k|k}\bigr]. 
 \label{Xmatrix:centered}
\end{align}

The main goal of this paper is to derive the {\it continuous-discrete} DF-EKF implementation methods. A traditional way for deriving the continuous-discrete filters is to apply numerical methods for solving the SDE in~\eqref{eq1.1}. This step yields the discretized stochastic system at hand and, next, one should derive the moment equations required by the filter. A standard routine for discretizing SDEs is to apply the Euler-Maruyama method, which is of strong order 0.5 as discussed in~\cite[Section~10.2]{KlPl99}. Some modern continuous-discrete Bayesian filters imply a higher order methods than the Euler-Maruyama integrator. For example, the It\^{o}-Taylor expansion of strong order 1.5 has been used to derive the continuous-discrete CKF in~\cite{Haykin2010}. In our work, we use both integrators and derive two {\it continuous-discrete} DF-EKF frameworks.

We start with the Euler-Maruyama method.  Let us consider the sampling interval $[t_{k-1}, t_{k}]$ and assume that it is additionally partitioned it into $L$ equally spaced subdivision nodes as follows: $t_{k-1}=t_{k-1}^{(0)} < t_{k-1}^{(1)} < \ldots t_{k-1}^{(l)}  < \ldots t_{k-1}^{(L)} = t_k$. In other words, the step size of the Euler-Maruyama integrator is $\delta=\Delta/L$ with $t_{k-1}^{(l)} = t_{k-1} + l \delta$, $l=0, \ldots, L-1$ and we have
\begin{align}
x_{k-1}^{(l+1)} & = f_d^{EM}\bigl(t_{k-1}^{(l)},x_{k-1}^{(l)}\bigr) + G \Delta w_l  \label{EM:scheme}
\end{align}
where the discretized drift function is
\begin{align}
f_d^{EM}\bigl(t_{k-1}^{(l)},x_{k-1}^{(l)}\bigr) & = x_{k-1}^{(l)} + \delta f\bigl(t_{k-1}^{(l)},x_{k-1}^{(l)}\bigr) \label{EM:scheme:drift}
\end{align}
and the random vector $\Delta w_l$ in~\eqref{EM:scheme} consists of the mutually independent and identically distributed Gaussian variables whose expected values are zero and their variances are $\delta Q$. More formally, they are generated from normal distribution ${\mathcal N}(\mathbf{0}, \delta Q)$.

The It\^{o}-Taylor expansion of strong order 1.5 applied to stochastic system~\eqref{eq1.1}, \eqref{eq1.2} on the mesh $t_{k-1}=t_{k-1}^{(0)} < t_{k-1}^{(1)} < \ldots t_{k-1}^{(l)}  < \ldots t_{k-1}^{(L)} = t_k$ over the sampling interval $[t_{k-1}, t_{k}]$ with step size  $\delta=\Delta/L$ has the following form~\cite[Section~10.4]{KlPl99}:
\begin{align}
x_{k-1}^{(l+1)} & = f_d^{IT}\bigl(t_{k-1}^{(l)},x_{k-1}^{(l)}\bigr) + GQ^{1/2} w_1 +{\mathbb L}f\bigl(t_{k-1}^{(l)},x_{k-1}^{(l)}\bigr) w_2 \label{IT:scheme}
\end{align}
where $Q^{1/2}$ is a square-root factor of the process covariance $Q = Q^{1/2}Q^{\top/2}$ and the discretized drift function is
{\small
\begin{align}
f_d^{IT}\bigl(t_{k-1}^{(l)},x_{k-1}^{(l)}\bigr) & = x_{k-1}^{(l)}+\delta f\bigl(t_{k-1}^{(l)},x_{k-1}^{(l)}\bigr) + \frac{1}{2}{\delta^2}{\mathbb L}_0 f\bigl(t_{k-1}^{(l)},x_{k-1}^{(l)}\bigr). \label{IT:scheme:drift}
\end{align}
}
The noises $w_1$ and $w_2$ in equation~\eqref{IT:scheme} are the correlated $n$-dimensional zero-mean random Gaussian variables with the properties discussed in~\cite{KlPl99}:
\begin{align}
w_1 & = \sqrt{\delta} v_1, \; \E{w_1w_1^\top} = \delta I, \; \E{w_1w_2^\top} = \frac{1}{2} \delta^2 I, \label{prop:1}\\
w_2 & = \frac{1}{2} \delta^{3/2}(v_1 + \frac{1}{\sqrt{3}}v_2), \; \E{w_2w_2^\top} =\frac{1}{3} \delta^3 I. \label{prop:2}
\end{align}

Finally, the term ${\mathbb L} f$ stands for a square matrix with $(i,j)$ entry being ${\mathbb L}_jf_i$. Let us denote the product $G^* = GQ^{1/2}$. The utilized differential operators ${\mathbb L}_0$ and ${\mathbb L}_j$ are defined as follows:
\begin{align}
{\mathbb L}_0 & = \frac{\partial}{\partial t} + \sum \limits_{i = 1}^{n} f_i \frac{\partial}{\partial x_i} + \frac{1}{2}\sum \limits_{j,p,r=1}^{n} G^*_{pj}  G^*_{rj} \frac{\partial^2}{\partial x_p \partial x_r}, \label{operator_L0} \\
{\mathbb L}_j &= \sum \limits_{i = 1}^{n} G^*_{ij} \frac{\partial}{\partial x_i}, \quad i,j = 1, 2,\ldots, n. \label{operator_Lj}
\end{align}

\begin{rmk}
It should be stressed that the numerical integrator based on the It\^{o}-Taylor expansion presented in formulas~\eqref{IT:scheme}~-- \eqref{operator_Lj} is valid for stochastic systems with time-invariant and state-independent diffusion term. For instance, if the diffusion $G(\cdot)$ related to the standard Brownian motion $\beta(t)$ in the examined state-space model or the covariance $Q(\cdot)$ are time-variant, then the differential operators ${\mathbb L}_0$ and ${\mathbb L}_j$ have a more sophisticated expression than formulas~\eqref{operator_L0}, \eqref{operator_Lj}, i.e. they do not hold in general. More details can be found in~\cite[Section~10.4]{KlPl99}.
\end{rmk}

\section{Derivation of the continuous-discrete derivative-free EKF implementation methods} \label{sec:main:conventional}

The goal of this section is to propose the continuous-discrete derivative-free EKF methods for stochastic system~\eqref{eq1.1}, \eqref{eq1.2}. We suppose that the system might be discretized by two numerical integrators presented in the previous section, i.e. either by the scheme in~\eqref{EM:scheme} or by equation~\eqref{IT:scheme}, respectively. Hence, we derive two filtering methods: 1) the Euler-Maruyama-based continuous-discrete DF-EKF (EM-0.5 DF-EKF), and 2) the continuous-discrete DF-EKF within the It\^{o}-Taylor expansion (IT-1.5 DF-EKF).

\subsection{Derivation of moment equations for prediction step}

We start with a derivation of the continuous-discrete EM-0.5 DF-EKF estimator.

\begin{lmm}\label{Lm1}
Let us explore the sampling interval $[t_{k-1}, t_{k}]$ with the equidistant mesh $t_{k-1}=t_{k-1}^{(0)} < t_{k-1}^{(1)} < \ldots t_{k-1}^{(l)}  < \ldots t_{k-1}^{(L)} = t_k$ where the integrator's step size is $\delta=\Delta/L$ and the nodes are $t_{k-1}^{(l)} = t_{k-1} + l \delta$, $l=0, \ldots, L-1$. Given the state estimate, $\hat x_{k-1|k-1}$, and the filter covariance matrix, $P_{k-1|k-1}$, at time $t_{k-1}$, set up the initial values for the integration scheme to be applied as follows: $\hat x_{k-1|k-1}^{(0)}:=\hat x_{k-1|k-1}$ and $P_{k-1|k-1}^{(0)}:=P_{k-1|k-1}$. If the Euler-Maruyama method in~\eqref{EM:scheme}, \eqref{EM:scheme:drift} is applied to the SDE in~\eqref{eq1.1}, then the derivative-free EKF framework yields the following state and covariance equations:
{\small
\begin{align}
\hat x_{k-1|k-1}^{(l+1)} & = f_d^{EM}\bigl(t_{k-1}^{(l)},\hat x_{k-1|k-1}^{(l)}\bigr)= \hat x_{k-1|k-1}^{(l)} + \delta f\bigl(t_{k-1}^{(l)}, \hat x_{k-1|k-1}^{(l)}), \label{EM:Predict:X} \\
P_{k-1|k-1}^{(l+1)} & = \overline{{\mathbb F}{\mathbb X}}_{k-1|k-1}^{EM,\;(l)}\left(\overline{{\mathbb F}{\mathbb X}}_{k-1|k-1}^{EM,\;(l)}\right)^{\top} + \delta GQG^{\top}  \label{EM:Predict:P}
\end{align}}
where the centered and scaled matrix $\overline{{\mathbb F}{\mathbb X}}_{k-1|k-1}^{EM,\;(l)}$ is defined from $n$ vectors ${\mathcal X}_{i, k-1|k-1}^{(l)}$, $i=1, \ldots,n$, generated around the current estimate, $\hat x_{k-1|k-1}^{(l)}$, by the rule
\begin{align}
{\mathbb  X}_{k-1|k-1}^{(l)} & = \hat x_{k-1|k-1}^{(l)} {\mathbf 1}^{\top} + \frac{\sqrt{n}}{\alpha} \left(P_{k-1|k-1}^{(l)}\right)^{1/2}, \label{need:proof1}
\end{align}
and propagated through the discretized drift function $f_d^{EM}(\cdot)$:
\begin{align}
{\mathbb F}{\mathbb X}_{k-1|k-1}^{EM,\;(l)} & =  f_d^{EM}\bigl(t_{k-1}^{(l)}, {\mathbb  X}_{k-1|k-1}^{(l)}\bigr),  \label{need:proof2}  \\
\overline{{\mathbb F}{\mathbb X}}_{k-1|k-1}^{EM,\;(l)} & = \frac{\alpha}{\sqrt{n}} \left[{\mathbb F}{\mathbb X}_{k-1|k-1}^{EM,\;(l)} - f_d^{EM}\bigl(t_{k-1}^{(l)},\hat x_{k-1|k-1}^{(l)}\bigr)\;{\mathbf 1}^{\top} \right].  \label{need:proof3}
\end{align}

\end{lmm}
\begin{pf}
Having applied the discretization scheme in~\eqref{EM:scheme} and taking into account the zero-mean random vector $\Delta w_l$, we derive
\begin{align}
\hat x_{k-1|k-1}^{(l+1)} & := \E{x_{k-1}^{(l+1)}} = \E{f_d^{EM}\bigl(t_{k-1}^{(l)},x_{k-1}^{(l)}\bigr) + G \Delta w_l}  \nonumber \\
& =\E{f_d^{EM}\bigl(t_{k-1}^{(l)},x_{k-1}^{(l)}\bigr)} = f_d^{EM}\bigl(t_{k-1}^{(l)},\hat x_{k-1|k-1}^{(l)} \bigr)\nonumber \\
& =\hat x_{k-1|k-1}^{(l)} + \delta f\bigl(t_{k-1}^{(l)}, \hat x_{k-1|k-1}^{(l)}), \label{EM:EKF:mean}
\end{align}
i.e. equation~\eqref{EM:Predict:X} for calculating the state estimate is proved.

Next, let us consider the difference
\begin{align}
x_{k-1}^{(l+1)} - \hat x_{k-1|k-1}^{(l+1)} & = f_d^{EM}\bigl(t_{k-1}^{(l)},x_{k-1}^{(l)}\bigr) + G \Delta w_l \nonumber \\
& - \E{f_d^{EM}\bigl(t_{k-1}^{(l)},x_{k-1}^{(l)}\bigr) + G \Delta w_l}. \label{proof1:new}
\end{align}

Taking into account~\eqref{proof1:new} and the fact that  $\Delta w_l$ is a zero-mean random vector with the covariance matrix $\delta Q$, we can calculate the covariance matrix as follows:
\begin{align}
P_{k-1|k-1}^{(l+1)} & = \E{\left[x_{k-1}^{(l+1)} - \hat x_{k-1|k-1}^{(l+1)} \right]\left[x_{k-1}^{(l+1)} - \hat x_{k-1|k-1}^{(l+1)}\right]^{\top}} \nonumber \\
& = {\mathbf E} \left\{ \left[f_d^{EM}\bigl(t_{k-1}^{(l)},x_{k-1}^{(l)}\bigr) - \E{f_d^{EM}\bigl(t_{k-1}^{(l)},x_{k-1}^{(l)}\bigr)} + G \Delta w_l  \right] \right.\nonumber \\
& \quad \times \left. \left[f_d^{EM}\bigl(t_{k-1}^{(l)},x_{k-1}^{(l)}\bigr) - \E{f_d^{EM}\bigl(t_{k-1}^{(l)},x_{k-1}^{(l)}\bigr)} + G \Delta w_l  \right]^{\top} \right\} \nonumber \\
& = \var{f_d^{EM}\bigl(t_{k-1}^{(l)},x_{k-1}^{(l)}\bigr)} + \delta G Q G^{\top} \label{trulll}
\end{align}
where the term $\var{f_d^{EM}\bigl(t_{k-1}^{(l)},x_{k-1}^{(l)}\bigr)}$ stands for the covariance matrix, i.e., we have
\begin{align}
& \var{f_d^{EM}\bigl(t_{k-1}^{(l)},x_{k-1}^{(l)}\bigr)}  = {\mathbf E}\left\{ \left[f_d^{EM}\bigl(t_{k-1}^{(l)},x_{k-1}^{(l)}\bigr) -  \E{f_d^{EM}\bigl(t_{k-1}^{(l)},x_{k-1}^{(l)}\bigr)}\right] \right.  \nonumber \\
& \quad \left. \times \left[f_d^{EM}\bigl(t_{k-1}^{(l)},x_{k-1}^{(l)}\bigr) - \E{f_d^{EM}\bigl(t_{k-1}^{(l)},x_{k-1}^{(l)}\bigr)}\right]^{\top}\right\} \nonumber \\
 & = {\mathbf E}\left\{ \left[f_d^{EM}\bigl(t_{k-1}^{(l)},x_{k-1}^{(l)}\bigr) -  f_d^{EM}\bigl(t_{k-1}^{(l)},\hat x_{k-1|k-1}^{(l)}\bigr)\right] \right.  \nonumber \\
 & \quad \left. \times \left[f_d^{EM}\bigl(t_{k-1}^{(l)},x_{k-1}^{(l)}\bigr) - f_d^{EM}\bigl(t_{k-1}^{(l)},\hat x_{k-1|k-1}^{(l)}\bigr)\right]^{\top}\right\}.
 \label{proof:main:new}
\end{align}

To compute the covariance matrix above, the derivative-free EKF principle is applied. More precisely, given the current values $\hat x_{k-1|k-1}^{(l)}$ and $P_{k-1|k-1}^{(l)}$, one generates a set of $n$ vectors $\xi_i$, $i = 1, \ldots, n$, that approximate the second moment by
\begin{align*}
P_{k-1|k-1}^{(l)}& = \frac{1}{n} \sum \limits_{i=1}^{n} \xi_i \xi_i^{\top},  & i & = 1, \ldots, n,
\end{align*}
and, next, a set of $n$ vectors ${\mathcal X}_{i, k-1|k-1}^{(l)}$, $i=1, \ldots,n$, are generated by the rule in~\eqref{Xi_vectors}, i.e.
\begin{align*}
{\mathcal X}_{i,k-1|k-1}^{(l)} & =  \hat x_{k-1|k-1}^{(l)} + \frac{1}{\alpha} \xi_i, & i & = 1, \ldots, n.
\end{align*}
Similarly to~\eqref{Xi_vectors:matrix}, the above equation can be written in a vector-matrix form as follows:
\begin{align*}
{\mathbb  X}_{k-1|k-1}^{(l)} & = \hat x_{k-1|k-1}^{(l)} {\mathbf 1}^{\top} + \frac{\sqrt{n}}{\alpha} \left(P_{k-1|k-1}^{(l)}\right)^{1/2},
\end{align*}
i.e., formula~\eqref{need:proof1} is proved.

Again, the matrix ${\mathbb  X}_{k-1|k-1}^{(l)}$ is collected from the vectors ${\mathcal X}_{i, k-1|k-1}^{(l)}$, $i=1, \ldots,n$, which are located by columns:
\begin{align}
{\mathbb  X}_{k-1|k-1}^{(l)} & =\Bigl[{\mathcal X}_{1,k-1|k-1}^{(l)} \;|\; \ldots \;|\; {\mathcal X}_{n,k-1|k-1}^{(l)}\Bigr]. \label{proof3:new}
\end{align}

Next, all vectors ${\mathcal X}_{i, k-1|k-1}^{(l)}$, $i=1, \ldots,n$, are propagated through the discretized drift function in~\eqref{EM:scheme}, i.e., we get the matrix collected from the propagated vectors as follows:
\begin{align}
& {\mathbb F}{\mathbb X}_{k-1|k-1}^{EM,\;(l)} = f_d^{EM}\bigl(t_{k-1}^{(l)}, {\mathbb  X}_{k-1|k-1}^{(l)}\bigr) \nonumber \\
& \quad= \bigl[f_d^{EM}\bigl(t_{k-1}^{(l)}, {\cal  X}_{1, k-1|k-1}^{(l)}\bigr) \;| \ldots |\; f_d^{EM}\bigl(t_{k-1}^{(l)}, {\cal  X}_{n, k-1|k-1}^{(l)}\bigr) \bigr],  \label{FX:matrix:center:1}
\end{align}
that is, exactly formula~\eqref{need:proof2}.

To calculate the last term in formula~\eqref{proof:main:new}, we also need to introduce the centered and scaled version of the matrix in~\eqref{FX:matrix:center:1}, i.e., we define
\begin{align}
\overline{{\mathbb F}{\mathbb X}}_{k-1|k-1}^{EM,\;(l)} & := \frac{\alpha}{\sqrt{n}} \left[f_d^{EM}\bigl(t_{k-1}^{(l)}, {\cal  X}_{1, k-1|k-1}^{(l)}\bigr)-f_d^{EM}\bigl(t_{k-1}^{(l)},\hat x_{k-1|k-1}^{(l)}) \;| \ldots \right. \nonumber \\
 & \left. \ldots |\; f_d^{EM}\bigl(t_{k-1}^{(l)}, {\cal  X}_{n, k-1|k-1}^{(l)}\bigr) - f_d^{EM}\bigl(t_{k-1}^{(l)},\hat x_{k-1|k-1}^{(l)}) \right] \nonumber \\
 & = \frac{\alpha}{\sqrt{n}} \left[f_d^{EM}\bigl(t_{k-1}^{(l)}, {\mathbb  X}_{k-1|k-1}^{(l)} \bigr)-f_d^{EM}\bigl(t_{k-1}^{(l)},\hat x_{k-1|k-1}^{(l)}){\mathbf 1}^{\top} \right] \nonumber \\
& = \frac{\alpha}{\sqrt{n}} \left[{\mathbb F}{\mathbb X}_{k-1|k-1}^{EM,\;(l)} - f_d^{EM}\bigl(t_{k-1}^{(l)},\hat x_{k-1|k-1}^{(l)}){\mathbf 1}^{\top} \right].
   \label{FX:matrix:center}
\end{align}
Thus, formula~\eqref{need:proof3} is also proved.

Following the derivative-free EKF principle for calculating the covariance matrix in~\eqref{P:approximation} and having substituted notations in~\eqref{FX:matrix:center}, we next obtain
\begin{align}
& \var{f_d^{EM}\bigl(t_{k-1}^{(l)},x_{k-1}^{(l)}\bigr)}   \nonumber \\
 & = {\mathbf E}\left\{ \left[f_d^{EM}\bigl(t_{k-1}^{(l)},x_{k-1}^{(l)}\bigr) -  f_d^{EM}\bigl(t_{k-1}^{(l)},\hat x_{k-1|k-1}^{(l)}\bigr)\right] \right.  \nonumber \\
 & \quad \left. \times \left[f_d^{EM}\bigl(t_{k-1}^{(l)},x_{k-1}^{(l)}\bigr) - f_d^{EM}\bigl(t_{k-1}^{(l)},\hat x_{k-1|k-1}^{(l)}\bigr)\right]^{\top}\right\} \nonumber \\
 & =  \frac{\alpha^2}{n} \sum \limits_{i=1}^{n}  \left\{ \left[f_d^{EM}\bigl(t_{k-1}^{(l)},  {\cal  X}_{i, k-1|k-1}^{(l)} \bigr) -  f_d^{EM}\bigl(t_{k-1}^{(l)},\hat x_{k-1|k-1}^{(l)}\bigr)\right] \right.  \nonumber \\
 & \quad \left. \times \left[f_d^{EM}\bigl(t_{k-1}^{(l)},  {\cal  X}_{i, k-1|k-1}^{(l)} \bigr) - f_d^{EM}\bigl(t_{k-1}^{(l)},\hat x_{k-1|k-1}^{(l)}\bigr)\right]^{\top}\right\} \nonumber \\
 & = \overline{{\mathbb F}{\mathbb X}}_{k-1|k-1}^{EM,\;(l)}\left(\overline{{\mathbb F}{\mathbb X}}_{k-1|k-1}^{EM,\;(l)}\right)^{\top}.
 \label{proof:main:new1}
\end{align}

Having substituted~\eqref{proof:main:new1} into equation~\eqref{trulll}, we arrive at
\begin{align*}
P_{k-1|k-1}^{(l+1)} & = \var{f_d^{EM}\bigl(t_{k-1}^{(l)},x_{k-1}^{(l)}\bigr)} + \delta G Q G^{\top} \nonumber \\
& = \overline{{\mathbb F}{\mathbb X}}_{k-1|k-1}^{EM,\;(l)}\left(\overline{{\mathbb F}{\mathbb X}}_{k-1|k-1}^{EM,\;(l)}\right)^{\top} + \delta G Q G^{\top},
\end{align*}
that is, exactly formula~\eqref{EM:Predict:P}.

This completes the proof of Lemma~\ref{Lm1}.
\end{pf}

Similarly, we next derive the prediction step of the continuous-discrete IT-1.5 DF-EKF estimator.

\begin{lmm}\label{Lm2}
Let us explore the sampling interval $[t_{k-1}, t_{k}]$ with the equidistant mesh $t_{k-1}=t_{k-1}^{(0)} < t_{k-1}^{(1)} < \ldots t_{k-1}^{(l)}  < \ldots t_{k-1}^{(L)} = t_k$ where the integrator's step size is $\delta=\Delta/L$ and the nodes are $t_{k-1}^{(l)} = t_{k-1} + l \delta$, $l=0, \ldots, L-1$. Given the state estimate, $\hat x_{k-1|k-1}$, and the filter covariance matrix, $P_{k-1|k-1}$, at time $t_{k-1}$, set up the initial values for the integration scheme to be applied as follows: $\hat x_{k-1|k-1}^{(0)}:=\hat x_{k-1|k-1}$ and $P_{k-1|k-1}^{(0)}:=P_{k-1|k-1}$. If the It\^{o}-Taylor expansion summarized in~\eqref{IT:scheme}~-- \eqref{operator_Lj} is applied to the SDE in~\eqref{eq1.1}, then the derivative-free EKF concept yields the following state and covariance equations:
\begin{align}
\hat x_{k-1|k-1}^{(l+1)} & = f_d^{IT}\bigl(t_{k-1}^{(l)},\hat x_{k-1|k-1}^{(l)}\bigr) = \hat x_{k-1|k-1}^{(l)} \nonumber \\
&  + \delta f\bigl(t_{k-1}^{(l)},\hat x_{k-1|k-1}^{(l)}\bigr) + \frac{1}{2}{\delta^2}{\mathbb L}_0 f\bigl(t_{k-1}^{(l)},\hat x_{k-1|k-1}^{(l)}\bigr), \label{IT:Predict:X} \\
P_{k-1|k-1}^{(l+1)} & = \overline{{\mathbb F}{\mathbb X}}_{k-1|k-1}^{IT,\;(l)}\left(\overline{{\mathbb F}{\mathbb X}}_{k-1|k-1}^{IT,\;(l)}\right)^{\top} + \delta GQG^{\top}\nonumber \\
& +\frac{\delta^2}{2}\biggl[GQ^{1/2} \;{\mathbb L} f^{\top}\!\!\bigl(t_{k-1}^{(l)},\hat x_{k-1|k-1}^{(l)}\bigr)  \nonumber \\
& + {\mathbb L}  f\bigl(t_{k-1}^{(l)},\hat x_{k-1|k-1}^{(l)}\bigr)Q^{\top/2}G^{\top}\biggr] \nonumber \\
& +\frac{\delta^3}{3}\!\bigl[{\mathbb L} f\bigl(t_{k-1}^{(l)},\hat x_{k-1|k-1}^{(l)}\bigr)\!\bigr]\bigl[{\mathbb L} f\bigl(t_{k-1}^{(l)},\hat x_{k-1|k-1}^{(l)}\bigr)\!\bigr]^{\top}  \label{IT:Predict:P}
\end{align}
where the operators ${\mathbb L}_0$  and ${\mathbb L}$ are calculated in~\eqref{operator_L0} and \eqref{operator_Lj}, respectively. The centered and scaled matrix $\overline{{\mathbb F}{\mathbb X}}_{k-1|k-1}^{IT,\;(l)}$ is defined from $n$ vectors ${\mathcal X}_{i, k-1|k-1}^{(l)}$, $i=1, \ldots,n$, generated around the current estimate, $\hat x_{k-1|k-1}^{(l)}$, by the rule
\begin{align}
{\mathbb  X}_{k-1|k-1}^{(l)} & = \hat x_{k-1|k-1}^{(l)} {\mathbf 1}^{\top} + \frac{\sqrt{n}}{\alpha} \left(P_{k-1|k-1}^{(l)}\right)^{1/2}, \label{ITneed:proof1}
\end{align}
and propagated through the discretized drift function $f_d^{IT}(\cdot)$:
\begin{align}
{\mathbb F}{\mathbb X}_{k-1|k-1}^{IT,\;(l)} & =  f_d^{IT}\bigl(t_{k-1}^{(l)}, {\mathbb  X}_{k-1|k-1}^{(l)}\bigr),  \label{ITneed:proof2}  \\
\overline{{\mathbb F}{\mathbb X}}_{k-1|k-1}^{IT,\;(l)} & = \frac{\alpha}{\sqrt{n}} \left[{\mathbb F}{\mathbb X}_{k-1|k-1}^{IT,\;(l)} - f_d^{IT}\bigl(t_{k-1}^{(l)},\hat x_{k-1|k-1}^{(l)}\bigr)\;{\mathbf 1}^{\top} \right].  \label{ITneed:proof3}
\end{align}

\end{lmm}
\begin{pf}
Having applied the discretization scheme in~\eqref{IT:scheme} and taking into account the zero-mean noises $w_1$ and $w_2$, we derive
\begin{align}
\hat x_{k-1|k-1}^{(l+1)} & := \E{x_{k-1}^{(l+1)}} = {\mathbf E} \bigl\{ f_d^{IT}\bigl(t_{k-1}^{(l)},x_{k-1}^{(l)}\bigr) + GQ^{1/2} w_1   \nonumber \\
& +{\mathbb L}f\bigl(t_{k-1}^{(l)},x_{k-1}^{(l)}\bigr) w_2 \bigr\} =\E{f_d^{IT}\bigl(t_{k-1}^{(l)},x_{k-1}^{(l)}\bigr)} \nonumber \\
& = f_d^{IT}\bigl(t_{k-1}^{(l)},\hat x_{k-1|k-1}^{(l)} \bigr) = \hat x_{k-1|k-1}^{(l)}  \nonumber \\
& + \delta f\bigl(t_{k-1}^{(l)},\hat x_{k-1|k-1}^{(l)}\bigr) + \frac{1}{2}{\delta^2}{\mathbb L}_0 f\bigl(t_{k-1}^{(l)},\hat x_{k-1|k-1}^{(l)}\bigr), \label{IT:EKF:mean}
\end{align}
i.e., equation~\eqref{IT:Predict:X} for calculating the state estimate is proved.

Let us consider the difference
\begin{align}
& x_{k-1}^{(l+1)} - \hat x_{k-1|k-1}^{(l+1)}  = f_d^{IT}\bigl(t_{k-1}^{(l)},x_{k-1}^{(l)}\bigr) + GQ^{1/2} w_1 + {\mathbb L}f\bigl(t_{k-1}^{(l)},x_{k-1}^{(l)}\bigr) w_2   \nonumber \\
& - \E{ f_d^{IT}\bigl(t_{k-1}^{(l)},x_{k-1}^{(l)}\bigr) + GQ^{1/2} w_1 + {\mathbb L}f\bigl(t_{k-1}^{(l)},x_{k-1}^{(l)}\bigr) w_2 }. \label{proof1:new:1}
\end{align}

Taking into account~\eqref{proof1:new:1} and the zero-mean noises $w_1$ and $w_2$ with the properties in~\eqref{prop:1}, \eqref{prop:2}, we can derive the formula for covariance matrix calculation as follows:
\begin{align}
P_{k-1|k-1}^{(l+1)} & = \E{\left[x_{k-1}^{(l+1)} - \hat x_{k-1|k-1}^{(l+1)} \right]\left[x_{k-1}^{(l+1)} - \hat x_{k-1|k-1}^{(l+1)}\right]^{\top}} \nonumber \\
& = \var{f_d^{IT}\bigl(t_{k-1}^{(l)},x_{k-1}^{(l)}\bigr)} + \delta G Q G^{\top} \nonumber \\
& +\frac{\delta^2}{2}\biggl[GQ^{1/2} \;{\mathbb L} f^{\top} \bigl(t_{k-1}^{(l)},\hat x_{k-1|k-1}^{(l)}\bigr)  \nonumber \\
& + {\mathbb L}  f\bigl(t_{k-1}^{(l)},\hat x_{k-1|k-1}^{(l)}\bigr)Q^{\top/2}G^{\top}\biggr] \nonumber \\
& +\frac{\delta^3}{3}\!\bigl[{\mathbb L} f\bigl(t_{k-1}^{(l)},\hat x_{k-1|k-1}^{(l)}\bigr)\!\bigr]\bigl[{\mathbb L} f\bigl(t_{k-1}^{(l)},\hat x_{k-1|k-1}^{(l)}\bigr)\!\bigr]^{\top} \label{trulll:1}
\end{align}
where the operator ${\mathbb L}$ is calculated by formula~\eqref{operator_Lj} and the term $\var{f_d^{IT}\bigl(t_{k-1}^{(l)},x_{k-1}^{(l)}\bigr)}$ stands for covariance matrix.

Similarly to equations~\eqref{proof:main:new}~-- \eqref{proof:main:new1} in the proof of Lemma~\ref{Lm1}, we show that
\begin{align}
\var{f_d^{IT}\bigl(t_{k-1}^{(l)},x_{k-1}^{(l)}\bigr)} & = \overline{{\mathbb F}{\mathbb X}}_{k-1|k-1}^{IT,\;(l)}\left(\overline{{\mathbb F}{\mathbb X}}_{k-1|k-1}^{IT,\;(l)}\right)^{\top}
 \label{proof:main:new1new}
\end{align}
where $\overline{{\mathbb F}{\mathbb X}}_{k-1|k-1}^{IT,\;(l)}$ stands for the scaled and centered matrix in~\eqref{ITneed:proof3} calculated by the DF-EKF concept in~\eqref{ITneed:proof1}, \eqref{ITneed:proof2}.

Having substituted result in~\eqref{proof:main:new1new} into equation~\eqref{trulll:1}, we prove formula~\eqref{IT:Predict:P}. This completes the proof of Lemma~\ref{Lm2}.
\end{pf}

\subsection{Measurement update step of the derivative-free EKF}

The measurement update step coincides with that proposed previously for discrete-time stochastic systems in~\cite{quine2006derivative}. More precisely, given $\hat x_{k|k-1}$ and $P_{k|k-1}$ calculated at the prediction step, one generates a set of $n$ vectors around the estimate $\hat x_{k|k-1}$ as follows:
\begin{align}
{\mathcal X}_{i,k|k-1} & =  \hat x_{k|k-1} + \frac{1}{\alpha} \xi_i, & i & = 1, \ldots, n, \label{ksi:new}
\end{align}
where the scalar parameter $\alpha$ is used to scale the discrete distribution about the mean.

The error covariance $P_{k|k-1}$ might be recovered from vectors ${\mathcal X}_{i, k|k-1}$, $i = 1, \ldots, n$, by taking into account the formula above:
\begin{align*}
P_{k|k-1} & = \frac{1}{n} \sum \limits_{i=1}^{n} \xi_i \xi_i^{\top} = \frac{\alpha^2}{n} \sum \limits_{i=1}^{n} [{\mathcal X}_{i, k|k-1} - \hat x_{k|k-1}][{\mathcal X}_{i, k|k-1} - \hat x_{k|k-1}]^{\top}.
\end{align*}

The vectors ${\mathcal X}_{i, k|k-1}$, $i=1, \ldots,n$  and the predicted state vector $\hat x_{k|k-1}$ are then propagated through the observation function $h(\cdot)$ to generate a set of vectors:
\begin{align}
Z_{i, k|k-1} & = h\left(k,{\mathcal X}_{i, k|k-1}\right), & i & = 1, \ldots, n, \label{eq:Zvectors} \\
\hat z_{k|k-1} & = h\left(k,\hat x_{k|k-1}\right). \label{eq:Zmean}
\end{align}

Next, the residual covariance and cross-covariance matrices are calculated by using vectors $\hat z_{k|k-1}$ and $Z_{i, k|k-1}$, $i = 1, \ldots, n$, as follows~\cite{quine2006derivative}:
\begin{align}
R_{e,k} & = \frac{\alpha^2}{n}\sum_{i=1}^{n} \bigl[Z_{i, k|k-1} -  \hat z_{k|k-1} \bigr] \bigl[Z_{i, k|k-1} -  \hat z_{k|k-1}\bigr]^\top + R_k, \label{eq:approx:Rek} \\
P_{xz,k} & = \frac{\alpha^2}{n}\sum_{i=1}^{n} \bigl[{\mathcal X}_{i, k|k-1} -  \hat x_{k|k-1} \bigr] \bigl[Z_{i, k|k-1} -  \hat z_{k|k-1}\bigr]^\top. \label{eq:approx:Pxz}
\end{align}

It is proved in~\cite{quine2006derivative} that in the limiting case as $\alpha \to \infty$, the derivative-free EKF formulas~\eqref{eq:approx:Rek}, \eqref{eq:approx:Pxz} converge to the standard EKF equations for calculating the  residual covariance $R_{e,k}$ and cross-covariance $P_{xz,k}$, respectively. The choice of the parameter $\alpha$ is adaptable and can be tuned on the application model at hand. More details can be found in the cited paper. The simulation results provided there show an excellent convergence to the standard EKF when $\alpha = 1000$.

Finally, the filtered estimate and the error covariance matrix of the derivative-free EKF are calculated as in the standard EKF, i.e. by the formulas (also see the algorithm in~\cite[Section~5]{quine2006derivative}):
\begin{align}
\hat x_{k|k} & =\hat x_{k|k-1}+{K}_k(z_k-\hat z_{k|k-1}),  \label{ckf:state}\\
{K}_{k} & = P_{xz,k}R_{e,k}^{-1}, \mbox{ and } P_{k|k}  = P_{k|k-1} - {K}_k R_{e,k} {K}_k^{\top}. \label{ckf:gain}
\end{align}

Again, we can show that formula~\eqref{ksi:new} can be written in a simple vector-matrix form
\begin{align}
{\mathbb  X}_{k|k-1}=\hat x_{k|k-1} {\mathbf 1}^{\top} + \frac{\sqrt{n}}{\alpha}P_{k|k-1}^{1/2} \label{newnewnewformula}
\end{align}
where ${\mathbf 1}$ is a column vector, whose elements are all equal to one and the matrix ${\mathbb  X}_{k|k-1}$ is collected from vectors ${\mathcal X}_{i, k|k-1}$, $i=1,\ldots, n$, by columns, i.e. we have
\begin{align*}
{\mathbb  X}_{k|k-1} & =\Bigl[{\mathcal X}_{1,k|k-1} \;|\; \ldots \;|\; {\mathcal X}_{n,k|k-1}\Bigr].
\end{align*}

Having introduced the centered and scaled matrices
 \begin{align*}
\overline{\mathbb  X}_{k|k-1} & = \frac{\alpha}{\sqrt{n}} \bigl[{\mathbb  X}_{k|k-1} -  \hat x_{k|k-1}{\mathbf 1}^{\top} \bigr] \\
& = \frac{\alpha}{\sqrt{n}} \bigl[{\mathcal X}_{1,k|k-1} - \hat x_{k|k-1} | \ldots | {\mathcal X}_{n,k|k-1} - \hat x_{k|k-1}\bigr], \\
\overline{\mathbb  Z}_{k|k-1} & =\frac{\alpha}{\sqrt{n}} \bigl[ {\mathbb  Z}_{k|k-1} -  \hat z_{k|k-1}{\mathbf 1}^{\top} \bigr] \\
& =\frac{\alpha}{\sqrt{n}}\bigl[{Z}_{1,k|k-1} - \hat z_{k|k-1}| \ldots | {Z}_{n,k|k-1} - \hat z_{k|k-1}\bigr].
\end{align*}
we write down equations~\eqref{eq:approx:Rek}, \eqref{eq:approx:Pxz} in a simple matrix form:
\begin{align}
R_{e,k} & =  \overline{\mathbb  Z}_{k|k-1} \overline{\mathbb  Z}_{k|k-1}^\top + R_k, & P_{xz,k} & = \overline{\mathbb  X}_{k|k-1} \overline{\mathbb  Z}_{k|k-1}^\top. \label{eq:approx:Pxz:new}
\end{align}

Finally, based on Lemma~\ref{Lm1} and~\ref{Lm2} with the measurement update formulas~\eqref{ckf:state}~-- \eqref{eq:approx:Pxz:new}, we formulate two novel continuous-discrete derivative-free EKF methods that are the EM-0.5 DF-EKF and IT-1.5 DF-EKF summarized in the left and right panels in Table~\ref{Tab:1}, respectively.

\begin{table*}[ht!]
{\small
\renewcommand{\arraystretch}{1.3}
\caption{The {\it conventional} continuous-discrete derivative-free EKF methods within the Euler-Maruyama discretization and It\^{o}-Taylor expansion.} \label{Tab:1}
\centering
\begin{tabular}{l||l|l}
\hline
& \cellcolor{myGray} {\bf Algorithm~1: EM-0.5 DF-EKF} & \cellcolor{myGray} {\bf Algorithm~2: IT-1.5 DF-EKF} \\
\hline
\hline
\textsc{Initialization:}  &  \multicolumn{2}{l}{0. Set the initial values $\hat x_{0|0} = \bar x_0$, $P_{0|0} = \Pi_0$ and derivative-free EKF parameter $\alpha$; e.g. $\alpha = 1000$ suggested in~\cite{quine2006derivative}.}\\
\hline
 \textsc{Time} & \multicolumn{2}{l}{1. On $[t_{k-1}, t_k]$ introduce a mesh: $t_{k-1}^{(l)} = t_{k-1} + l \delta$, $l=0, \ldots, L-1$,  $\delta=\Delta/L$, $\Delta = |t_k-t_{k-1}|$, i.e. $t_{k-1}=t_{k-1}^{(0)} < \ldots t_{k-1}^{(l)}  < \ldots t_{k-1}^{(L)} = t_k$.} \\
\textsc{Update (TU):} & \multicolumn{2}{l}{\; Set up the initial values for integrators: $\hat x_{k-1|k-1}^{(0)}:=\hat x_{k-1|k-1}$ and $P_{k-1|k-1}^{(0)}:=P_{k-1|k-1}$ at time node $t_{k-1}$.} \\
& \multicolumn{2}{l}{\; {\bf For} $l=1,\ldots L-1$ perform the following steps:}\\
& \multicolumn{2}{l}{\quad 2. Cholesky decomposition: $P_{k-1|k-1}^{(l)}=\left(P_{k-1|k-1}^{(l)}\right)^{1/2}\left(P_{k-1|k-1}^{(l)}\right)^{\top/2}$. Define ${\mathbb  X}_{k-1|k-1}^{(l)} = \hat x_{k-1|k-1}^{(l)} {\mathbf 1}^{\top} + \frac{\sqrt{n}}{\alpha} \left(P_{k-1|k-1}^{(l)}\right)^{1/2}$.} \\
& \quad 3. Calculate estimate $\hat x_{k-1|k-1}^{(l+1)}$ by eq.~\eqref{EM:Predict:X}. & \quad 3. Calculate estimate $\hat x_{k-1|k-1}^{(l+1)}$ by eq.~\eqref{IT:Predict:X}. \\
 & \quad 4. Propagate sample points: ${\mathbb F}{\mathbb X}_{k-1|k-1}^{EM,\;(l)} =  f_d^{EM}\bigl(t_{k-1}^{(l)}, {\mathbb  X}_{k-1|k-1}^{(l)}\bigr)$. &
 \quad 4. Propagate ${\mathbb F}{\mathbb X}_{k-1|k-1}^{IT,\;(l)} =  f_d^{IT}\bigl(t_{k-1}^{(l)}, {\mathbb  X}_{k-1|k-1}^{(l)}\bigr)$.  \\
& \quad 5. Find $\overline{{\mathbb F}{\mathbb X}}_{k-1|k-1}^{EM,\;(l)} = \frac{\alpha}{\sqrt{n}} \left[{\mathbb F}{\mathbb X}_{k-1|k-1}^{EM,\;(l)} - f_d^{EM}\bigl(t_{k-1}^{(l)},\hat x_{k-1|k-1}^{(l)}\bigr)\;{\mathbf 1}^{\top} \right]$. & \quad 5. $\overline{{\mathbb F}{\mathbb X}}_{k-1|k-1}^{IT,\;(l)} = \frac{\alpha}{\sqrt{n}} \left[{\mathbb F}{\mathbb X}_{k-1|k-1}^{IT,\;(l)} - f_d^{IT}\bigl(t_{k-1}^{(l)},\hat x_{k-1|k-1}^{(l)}\bigr)\;{\mathbf 1}^{\top} \right]$.  \\
& \quad 6. Calculate covariance $P_{k-1|k-1}^{(l+1)}$ by eq.~\eqref{EM:Predict:P}. & \quad 6. Calculate covariance $P_{k-1|k-1}^{(l+1)}$ by eq.~\eqref{IT:Predict:P}. \\
& \multicolumn{2}{l}{\; {\bf End}. At the last node, set $\hat x_{k|k-1}:=\hat x_{k-1|k-1}^{(L)}$ and $P_{k|k-1}:=P_{k-1|k-1}^{(L)}$ at time node $t_{k-1}^{(L)}:=t_k$.}\\
\hline
\textsc{Measurement} & \multicolumn{2}{l}{7. Perform Cholesky decomposition $P_{k|k-1}=P_{k|k-1}^{1/2}P_{k|k-1}^{\top/2}$.  Define sample points ${\mathbb  X}_{k|k-1}=\hat x_{k|k-1} {\mathbf 1}^{\top} + \frac{\sqrt{n}}{\alpha}P_{k|k-1}^{1/2}$.} \\
\textsc{Update (MU):} & \multicolumn{2}{l}{8. Propagate $\hat z_{k|k-1} = h\left(k,\hat x_{k|k-1}\right)$ and all sample points ${\mathbb  Z}_{k|k-1} = h\left(k,{\mathbb  X}_{k|k-1}\right)$.} \\
& \multicolumn{2}{l}{9.  Define the scaled and centered matrices: $\overline{\mathbb  X}_{k|k-1} = \frac{\alpha}{\sqrt{n}} \bigl[{\mathbb  X}_{k|k-1} -  \hat x_{k|k-1}{\mathbf 1}^{\top} \bigr]$ and $\overline{\mathbb  Z}_{k|k-1}  =\frac{\alpha}{\sqrt{n}} \bigl[ {\mathbb  Z}_{k|k-1} -  \hat z_{k|k-1}{\mathbf 1}^{\top} \bigr]$.  } \\
& \multicolumn{2}{l}{10. Find $R_{e,k}=\overline{\mathbb Z}_{k|k-1}\overline{\mathbb Z}_{k|k-1}^{\top}+R_k$, $P_{xz,k}=\overline{\mathbb X}_{k|k-1}\overline{\mathbb Z}_{k|k-1}^{\top}$, ${K}_{k}=P_{xz,k}R_{e,k}^{-1}$, $\hat x_{k|k}=\hat x_{k|k-1}+{K}_k(z_k-\hat z_{k|k-1})$, $P_{k|k}=P_{k|k-1} - {K}_k R_{e,k} {K}_k^{\top}$.} \\
\hline
\end{tabular}
}
\end{table*}

\begin{rmk}
The thoughtful readers may note that the IT-1.5 DF-EKF requires some derivatives calculation at prediction step inherent from the It\^{o}-Taylor {\it discretization} scheme; see equations~\eqref{operator_L0}, \eqref{operator_Lj}. However, we name the resulted filter as the IT-1.5 DF-EKF honouring the derivative-free principal utilized for covariance matrix approximation in such filters.
\end{rmk}

We next discuss some properties of Algorithms~1 and 2 proposed. Firstly, we note that the IT-1.5 DF-EKF method is developed within a higher order integration scheme than the EM-0.5 DF-EKF algorithm. This means that the IT-1.5 DF-EKF requires less subdivisions $L$ at the prediction step while integration process than the EM-0.5 DF-EKF to maintain a high estimation quality. The IT-1.5 DF-EKF is also expected to outperform the EM-0.5 DF-EKF for estimation accuracy.   Secondly, both methods derived are the {\it covariance}-type implementations because the error covariance matrix is processed while the filtering process. Recall, the information matrix, which is the inverse of the filter covariance, is predicted and updated in the {\it information}-type implementations. Thirdly, we note that the square-root factors of the filter covariance matrix are required at the prediction and measurement update steps to generate the sample points. We utilized the Cholesky decomposition in lines~2 and~7 to find the square-root factors, which is a traditional way to perform such operation in the realm of nonlinear Bayesian filters' implementation methods. The readers may refer to modern sample-data KF-like estimators, e.g., the CKF in~\cite{Haykin2009,Haykin2010,Haykin2018}, the UKF in~\cite{Julier2004,KnLe19}, the GHQF in~\cite{Ito2000} and many others. Unfortunately, all such methods suffer from numerical instability in a finite precision arithmetics; see also the discussion in~\cite[p.~1262]{Haykin2009}. Indeed, the roundoff errors may destroy the theoretical properties of the filter covariance matrix and make the Cholesky decomposition unfeasible, i.e. the filtering process is unexpectedly interrupted. Surely, the square-root factors might be found in an alternative way, for example, by using singular value decomposition (SVD), but the resulted conventional implementation methods are still vulnerable to roundoff; e.g., see the results of numerical tests in~\cite{KuKu20Automatica}. This means that they fail while practical computations due to roundoff and/or numerical integration errors. The best way to resolve this problem, i.e. to improve the numerical stability of the filtering methods, is to avoid the square-root operation entirely, which is possible within the square-root implementation way. In the next section, we derive such numerically stable algorithms for the EM-0.5 DF-EKF and IT-1.5 DF-EKF estimators derived in this paper.

\section{Derivation of square-root implementation methods for continuous-discrete derivative-free EKF framework} \label{sec:main:squareroot}

To avoid the square-root operation at each filtering step of the sample-data KF-like nonlinear estimators, one may follow the square-root strategy. The goal is to process the square-root factors of the filter covariance matrices $P_{k|k-1}^{1/2}$ and $P_{k|k}^{1/2}$ instead of the full matrices $P_{k|k-1}$ and $P_{k|k}$ in Algorithms~1 and~2. This means that the formulas of the novel EM-0.5 DF-EKF and IT-1.5 DF-EKF should be re-derived in terms of square-root factors in the following way: the matrix $\Pi_0 > 0$ is decomposed at the initial step and then the filtering methods propagate and update square-root factors of the covariance matrices, only. This routine ensures the positive (semi-) definiteness and symmetric form of the resulted error covariance $P_{k|k}^{1/2}P_{k|k}^{\top/2} = P_{k|k}$ at any time instance $t_k$ and yields the improved numerical stability to roundoff errors; see more detail in~\cite[Chapter~7]{2015:Grewal:book}. To start a derivation of any square-root implementation method, one needs to take into account the following circumstance: the factorization $\Pi_0 = \Pi_0^{1/2}\Pi_0^{\top/2}$ might be performed in various ways and, hence, a wide variety of the factored-form implementation methods comes from the chosen factorization strategy.

\subsection{Cholesky-based square-root implementation methods}

The traditional way of deriving the square-root implementation methods in the KF realm is based on an utilization of the Cholesky decomposition. Throughout the paper, we consider the lower triangular Cholesky factors, i.e. $P = SS^{\top}$ where $S$ is a lower triangular matrix with positive diagonal entries~\cite[Section~1.3.2]{Bj15}. In this section, we design the Cholesky factored-form derivative-free EKF methods with the lower triangular square-root matrix $S:=P_{k|k}^{1/2}$ at the measurement update step as well as the lower triangular factor $S:=P_{k|k-1}^{1/2}$ at the prediction filtering step. We also take into account that the utilization of orthogonal transformations for updating the square-root factors additionally improves the numerical stability to roundoff errors. In other words, it is desirable to use orthogonal rotations for updating the filters' quantities as far as possible. Recall, the orthogonal transformations are used for computing the Cholesky factorization of a positive definite matrix that obeys the equation $C = AA^{\top} + BB^{\top}$ by applying the orthogonal transformation to the pre-array $[A \; B]$ as follows: $[A \; B] Q = [R \; 0]$ where $R$ is a lower triangular Cholesky factor of the matrix $C$ as we are looking for. It can be easily proved by multiplying the pre- and post-arrays involved, i.e. $[A \; B] Q Q^{\top}[A \; B]^{\top} = RR^{\top}$, and comparing both sides of the obtained formulas.

Recently, the Cholesky-based implementation methods have been suggested for the derivative-free EKF methods within the MATLAB ODE's solves in~\cite{kulikova2023derivative}. The measurement updates of those implementation methods coincide with the related square-root EM-0.5 DF-EKF and IT-1.5 DF-EKF steps to be derived here. Thus, we briefly summarize this part and, next, derive the prediction step equations.

The first equation in~\eqref{eq:approx:Pxz:new} can be factorized as follows:
\[
R_{e,k} = \overline{\mathbb Z}_{k|k-1}\overline{\mathbb Z}_{k|k-1}^{\top} + R_k = \left[\overline{\mathbb Z}_{k|k-1} \;\; R_k^{1/2}\right]\left[\overline{\mathbb Z}_{k|k-1} \;\; R_k^{1/2}\right]^{\top}
\]
and, hence, we have the computational way for finding the square-root factor $R_{e,k}^{1/2}$ of the residual covariance by
\begin{equation}
\underbrace{\left[\overline{\mathbb Z}_{k|k-1} \quad R_k^{1/2}\right]}_{pre-array} \Theta = \underbrace{[R_{e,k}^{1/2} \quad 0]}_{post-array}
\label{SR:Rek:1}
\end{equation}
where $\Theta$ is any orthogonal transformation that lower triangulates the pre-array.

The formula for the filter's gain computation in terms of the square-root factor $R_{e,k}^{1/2}$ is as follows:
\[ K_{k} = P_{xz,k}R_{e,k}^{-1} = P_{xz,k}R_{e,k}^{-\top/2}R_{e,k}^{-1/2}. \]

The second equation in~\eqref{ckf:gain} for calculating the filter covariance $P_{k|k}$ can not be factorized straightforward, but the symmetric equation for calculating $P_{k|k}$ has been derived in~\cite{kulikova2023derivative}:
\[
P_{k|k} = \left[\overline{\mathbb X}_{k|k-1}-{K}_{k}\overline{\mathbb Z}_{k|k-1}\right]\left[\overline{\mathbb X}_{k|k-1}-{K}_{k}\overline{\mathbb Z}_{k|k-1}\right]^{\top} + K_k R_k K_k^{\top}.
\]
and, hence, the square-root factor $P_{k|k}^{1/2}$ can be found through
\begin{equation}
\underbrace{\left[ (\overline{\mathbb X}_{k|k-1}-{K}_{k}\overline{\mathbb Z}_{k|k-1}) \quad K_kR_k^{1/2}\right] }_{pre-array} \Theta = \underbrace{ [P_{k|k}^{1/2} \quad 0]}_{post-array} \label{SR:P:1}
 \end{equation}
where $\Theta$ is any orthogonal matrix that lower triangulates the pre-array.

Alternative formula for simultaneous calculation of the square-root factors $R_{e,k}^{1/2}$ and $P_{k|k}^{1/2}$ is given as follows\footnote{It can be easily proved by multiplying the pre- and post-arrays involved, and comparing both sides of the obtained formulas with equations~\eqref{ckf:gain} and~\eqref{eq:approx:Pxz:new}.}:
\begin{equation}
\underbrace{\begin{bmatrix}
\overline{\mathbb  Z}_{k|k-1}  & R_k^{1/2} \\
\overline{\mathbb  X}_{k|k-1}  & 0
\end{bmatrix}}_{pre-array} \Theta =
\underbrace{\begin{bmatrix}
R_{e,k}^{1/2} & 0\\
\bar P_{xz,k} & P^{1/2}_{k|k}
\end{bmatrix}}_{post-array} \label{SR:array}
\end{equation}
where $\Theta$ is any orthogonal matrix that (block) lower triangulates the pre-array, i.e. $R_{e,k}^{1/2}$ and $P_{k|k}^{1/2}$ are the lower triangular Cholesky factors, which we are looking for. They are simply read-off from the post-array.

Finally, we need to derive square-root formulas for the prediction step of the novel EM-0.5 DF-EKF and IT-1.5 DF-EKF methods in Algorithms~1 and~2. In other words, we need to factorize formulas~\eqref{EM:Predict:P} and~\eqref{IT:Predict:P}, which are proved in Lemmas~1 and~2, respectively. It is easy to show that equation~\eqref{EM:Predict:P} of the EM-0.5 DF-EKF method can be written in the form:
\begin{align*}
  P_{k-1|k-1}^{(l+1)} & = \overline{{\mathbb F}{\mathbb X}}_{k-1|k-1}^{EM,\;(l)}\left(\overline{{\mathbb F}{\mathbb X}}_{k-1|k-1}^{EM,\;(l)}\right)^{\top} + \delta GQG^{\top}\nonumber \\
  &  = \left[\overline{{\mathbb F}{\mathbb X}}_{k-1|k-1}^{EM,\;(l)} \quad \sqrt{\delta} GQ^{1/2}\right]
  \left[ \overline{{\mathbb F}{\mathbb X}}_{k-1|k-1}^{EM,\;(l)} \quad \sqrt{\delta} GQ^{1/2} \right]^{\top}
\end{align*}
where $Q^{1/2}$ is the lower triangular Cholesky factor of the process covariance matrix $Q$, i.e. $Q=Q^{1/2}Q^{\top/2}$.

Thus, the computational way for finding the square-root factor $\left(P_{k-1|k-1}^{(l+1)}\right)^{1/2}$ of the predicted filter covariance matrix is
\begin{equation}
\underbrace{\biggl[\overline{{\mathbb F}{\mathbb X}}_{k-1|k-1}^{EM,\;(l)} \quad \; \sqrt{\delta} GQ^{1/2} \biggr]}_{pre-array} \Theta = \underbrace{\biggl[\left(P_{k-1|k-1}^{(l+1)}\right)^{1/2} \quad 0\biggr]}_{post-array} \label{eq:EM:SR:P}
\end{equation}
where $\Theta$ is any orthogonal transformation that lower triangulates the pre-array. The equality above can be proved by multiplying the pre- and post-arrays involved and comparing both sides with formula~\eqref{EM:Predict:P}.

Similarly, we re-derive formula~\eqref{IT:Predict:P} of the IT-1.5 DF-EKF in terms of the square-root factors of the predicted filter covariance matrix. For that, we re-arrange equation~\eqref{IT:Predict:P}  as follows:
\begin{align}
& P_{k-1|k-1}^{(l+1)}  = \left[ \overline{{\mathbb F}{\mathbb X}}_{k-1|k-1}^{IT,\;(l)}, \; \sqrt{\delta}(GQ^{1/2}+\frac{\delta}{2}{\mathbb L}f_{k-1}^{(l)}), \; \sqrt{\frac{\delta^3}{12}} {\mathbb L}f_{k-1}^{(l)} \right] \nonumber \\
& \times \left[ \overline{{\mathbb F}{\mathbb X}}_{k-1|k-1}^{IT,\;(l)}, \; \sqrt{\delta}(GQ^{1/2}+\frac{\delta}{2}{\mathbb L}f_{k-1}^{(l)}), \; \sqrt{\frac{\delta^3}{12}}
{\mathbb L}f_{k-1}^{(l)}\right]^{\top} \label{eq:proof:1}
\end{align}
where ${\mathbb L}f_{k-1}^{(l)} := {\mathbb L} f\bigl(t_{k-1}^{(l)},\hat x_{k-1|k-1}^{(l)}\bigr)$. Hence, the square-root factor of the predicted filter covariance $P_{k-1|k-1}^{(l+1)}$ can be found by
\begin{align}
& \underbrace{\left[ \overline{{\mathbb F}{\mathbb X}}_{k-1|k-1}^{IT,\;(l)}, \; \sqrt{\delta}(GQ^{1/2}+\frac{\delta}{2}{\mathbb L}f_{k-1}^{(l)}), \; \sqrt{\frac{\delta^3}{12}} {\mathbb L}f_{k-1}^{(l)} \right] }_{pre-array} \Theta \nonumber \\
& = \underbrace{ \biggl[\left(P_{k-1|k-1}^{(l+1)}\right)^{1/2} \quad 0\biggr]}_{post-array} \label{eq:IT:SR:P}
\end{align}
where $\Theta$ is any orthogonal matrix that lower triangulates the pre-array. The equality above can be proved by multiplying the pre- and post-arrays involved and comparing both sides with formula~\eqref{IT:Predict:P}.

\begin{table*}[ht!]
{\small
\renewcommand{\arraystretch}{1.3}
\caption{The {\it Cholesky-based} continuous-discrete derivative-free EKF methods within the Euler-Maruyama discretization and It\^{o}-Taylor expansion.} \label{Tab:2}
\centering
\begin{tabular}{l||l|l}
\hline
& \cellcolor{myGray} {\bf Algorithm~1a: Cholesky-based EM-0.5 DF-EKF} & \cellcolor{myGray} {\bf Algorithm~2a: Cholesky-based IT-1.5 DF-EKF} \\
\hline
\hline
\textsc{Initialization:}  &  \multicolumn{2}{l}{0. Cholesky decomposition $\Pi_0=\Pi_0^{1/2}\Pi_0^{\top/2}$. Set $\hat x_{0|0} = \bar x_0$, $P_{0|0}^{1/2} = \Pi_0^{1/2}$ and parameter $\alpha$; e.g. $\alpha = 1000$ suggested in~\cite{quine2006derivative}.}\\
\hline
 \textsc{Time} & \multicolumn{2}{l}{1. On $[t_{k-1}, t_k]$ introduce a mesh: $t_{k-1}^{(l)} = t_{k-1} + l \delta$, $l=0, \ldots, L-1$,  $\delta=\Delta/L$, $\Delta = |t_k-t_{k-1}|$, i.e. $t_{k-1}=t_{k-1}^{(0)} < \ldots t_{k-1}^{(l)}  < \ldots t_{k-1}^{(L)} = t_k$.} \\
\textsc{Update (TU):} & \multicolumn{2}{l}{\; Set up the initial values for integrators: $\hat x_{k-1|k-1}^{(0)}:=\hat x_{k-1|k-1}$ and $\left(P_{k-1|k-1}^{(0)}\right)^{1/2}:=P_{k-1|k-1}^{1/2}$ at time node $t_{k-1}$.} \\
& \multicolumn{2}{l}{\; {\bf For} $l=1,\ldots L-1$ perform the following steps:}\\
& \multicolumn{2}{l}{\quad 2. Generate all sample points through ${\mathbb  X}_{k-1|k-1}^{(l)} = \hat x_{k-1|k-1}^{(l)} {\mathbf 1}^{\top} + \frac{\sqrt{n}}{\alpha} \left(P_{k-1|k-1}^{(l)}\right)^{1/2}$.} \\
& \quad 3. Calculate estimate $\hat x_{k-1|k-1}^{(l+1)}$ by eq.~\eqref{EM:Predict:X}. & \quad 3. Calculate estimate $\hat x_{k-1|k-1}^{(l+1)}$ by eq.~\eqref{IT:Predict:X}. \\
 & \quad 4. Propagate sample points: ${\mathbb F}{\mathbb X}_{k-1|k-1}^{EM,\;(l)} =  f_d^{EM}\bigl(t_{k-1}^{(l)}, {\mathbb  X}_{k-1|k-1}^{(l)}\bigr)$. &
 \quad 4. Propagate ${\mathbb F}{\mathbb X}_{k-1|k-1}^{IT,\;(l)} =  f_d^{IT}\bigl(t_{k-1}^{(l)}, {\mathbb  X}_{k-1|k-1}^{(l)}\bigr)$.  \\
& \quad 5. Find $\overline{{\mathbb F}{\mathbb X}}_{k-1|k-1}^{EM,\;(l)} = \frac{\alpha}{\sqrt{n}} \left[{\mathbb F}{\mathbb X}_{k-1|k-1}^{EM,\;(l)} - f_d^{EM}\bigl(t_{k-1}^{(l)},\hat x_{k-1|k-1}^{(l)}\bigr)\;{\mathbf 1}^{\top} \right]$. & \quad 5. $\overline{{\mathbb F}{\mathbb X}}_{k-1|k-1}^{IT,\;(l)} = \frac{\alpha}{\sqrt{n}} \left[{\mathbb F}{\mathbb X}_{k-1|k-1}^{IT,\;(l)} - f_d^{IT}\bigl(t_{k-1}^{(l)},\hat x_{k-1|k-1}^{(l)}\bigr)\;{\mathbf 1}^{\top} \right]$.  \\
& \quad 6. Calculate square-root $\left(P_{k-1|k-1}^{(l+1)}\right)^{1/2}$ by eq.~\eqref{eq:EM:SR:P}. & \quad 6. Calculate square-root $\left(P_{k-1|k-1}^{(l+1)}\right)^{1/2}$ by eq.~\eqref{eq:IT:SR:P}. \\
& \multicolumn{2}{l}{\; {\bf End}. At the last node, set $\hat x_{k|k-1}:=\hat x_{k-1|k-1}^{(L)}$ and $P_{k|k-1}^{1/2}:=\left(P_{k-1|k-1}^{(L)}\right)^{1/2}$ at time node $t_{k-1}^{(L)}:=t_k$.}\\
\hline
\textsc{Measurement} & \multicolumn{2}{l}{7. Define sample points ${\mathbb  X}_{k|k-1}=\hat x_{k|k-1} {\mathbf 1}^{\top} + \frac{\sqrt{n}}{\alpha}P_{k|k-1}^{1/2}$. Propagate $\hat z_{k|k-1} = h\left(k,\hat x_{k|k-1}\right)$ and ${\mathbb  Z}_{k|k-1} = h\left(k,{\mathbb  X}_{k|k-1}\right)$.} \\
\textsc{Update (MU):} & \multicolumn{2}{l}{8.  Define the scaled and centered matrices: $\overline{\mathbb  X}_{k|k-1} = \frac{\alpha}{\sqrt{n}} \bigl[{\mathbb  X}_{k|k-1} -  \hat x_{k|k-1}{\mathbf 1}^{\top} \bigr]$ and $\overline{\mathbb  Z}_{k|k-1}  =\frac{\alpha}{\sqrt{n}} \bigl[ {\mathbb  Z}_{k|k-1} -  \hat z_{k|k-1}{\mathbf 1}^{\top} \bigr]$.  } \\
& \multicolumn{2}{l}{9. Find $R_{e,k}^{1/2}$ by eq.~\eqref{SR:Rek:1}, $P_{xz,k}=\overline{\mathbb X}_{k|k-1}\overline{\mathbb Z}_{k|k-1}^{\top}$, ${K}_{k}=P_{xz,k}R_{e,k}^{-\top/2}R_{e,k}^{-1/2}$, $\hat x_{k|k}=\hat x_{k|k-1}+{K}_k(z_k-\hat z_{k|k-1})$ and $P_{k|k}^{1/2}$ by eq.~\eqref{SR:P:1}.} \\
\hline
\hline
& \cellcolor{myGray} {\bf Algorithm~1b: Cholesky-based EM-0.5 DF-EKF} & \cellcolor{myGray} {\bf Algorithm~2b: Cholesky-based IT-1.5 DF-EKF} \\
\textsc{Initialization:}  &  \qquad $\rightarrow$ Repeat from Algorithm~1a. &  \qquad $\rightarrow$ Repeat from Algorithm~2a.  \\
\textsc{Time Update} & \qquad $\rightarrow$ Repeat from Algorithm~1a. &  \qquad $\rightarrow$ Repeat from Algorithm~2a.  \\
\hline
\textsc{Measurement} &   \qquad $\rightarrow$ Repeat lines~7,8 of Algorithm~1a. &  \qquad $\rightarrow$ Repeat lines~7,8 of Algorithm~2a.  \\
\textsc{Update (MU):}  & \multicolumn{2}{l}{9. Apply transformation in~\eqref{SR:array} and read-off $R_{e,k}^{1/2}$, $\bar P_{xz,k}$ and $P^{1/2}_{k|k}$. Find
$K_k = \bar P_{xz,k} R_{e,k}^{-1/2}$ and $\hat x_{k|k}=\hat x_{k|k-1}+{K}_k(z_k-\hat z_{k|k-1})$.} \\
\hline
\end{tabular}
}
\end{table*}

To summarize, we suggest two variants of the Cholesky-based implementation methods for the EM-0.5 DF-EKF and IT-1.5 DF-EKF, respectively. They are presented in the form of pseudo-codes in Table~\ref{Tab:2} in Algorithms~1a, 1b and Algorithms~2a, 2b, respectively. As can be seen, Algorithms~1a and~2a require two QR factorizations at the measurement update steps. Meanwhile, Algorithms~1b and~2b demand only one QR transformation at each iterate, i.e. they are a bit faster than their counterparts in Algorithms~1a and~2a.

\subsection{SVD-based square-root implementation methods}

An alternative approach for designing the square-root filters is to apply the singular value decomposition. Recall, the rank-$r$ decomposition of matrix $A \in {\mathbb R}^{m\times n}$, $m<n$, where $r<m$ is given by~\cite[Theorem~1.1.6]{Bj15}:
$ A  = W\Sigma V^{\top}, \,
\Sigma =
\begin{bmatrix}
S & 0
\end{bmatrix} \in {\mathbb R}^{m\times n},  \; S={\rm diag}\{ \sigma_1,\ldots,\sigma_r\}
$ where $W \in {\mathbb R}^{m\times m}$, $V \in {\mathbb R}^{n\times n}$ are orthogonal matrices, and $\sigma_1\geq \ldots \geq\sigma_r>0$ are
the singular values of $A$.

In general, the SVD-based square-root factor $S$ of a symmetric matrix $P$ is defined by $S = W  \Sigma^{1/2}$ because of the decomposition $P=W  \Sigma W^{\top}$. Indeed, having applied SVD to the symmetric and full rank covariance matrix $\Pi_0$, we obtain $\Pi_0 = W_{\Pi_0} \Sigma_{\Pi_0} W_{\Pi_0}^{\top}$ where $W_{\Pi_0}$ is the orthogonal factor of size $n$, and $\Sigma_{\Pi_0}$ is a diagonal matrix of size $n$. Thus, the SVD-based square-root factor at the initial filtering step is given by $\Pi_0^{1/2} = W_{\Pi_0} \Sigma_{\Pi_0}^{1/2}$. Next, the novel EM-0.5 DF-EKF and IT-1.5 DF-EKF methods should be re-derived in terms of propagating and updating the orthogonal factors $W_{P_{k|k-1}}$, $W_{P_{k|k}}$ and diagonal factors $\Sigma_{P_{k|k-1}}$, $\Sigma_{P_{k|k}}$, respectively. For that, the SVD factorization is naturally used at the prediction and filtering steps of such filters' implementation methods.

It is also worth noting here that the SVD-based square-root filtering yields a square-root factor, which is a full square matrix (in general, it
might be a rectangular matrix), compared to a triangular factor obtained under the Cholesky-based decomposition approach discussed in previous section. The
SVD is also known to be the most accurate matrix factorization method, especially when the matrix to be factorized is close
to a singular one. Besides, SVD exists for any matrix that is not the case for Cholesky decomposition. Another benefit of any SVD-based square-root filtering algorithm is that all eigenvalues of the predicted and filtered error covariance matrices are available while estimation and this information might be used for automatic model analysis (e.g. detection of singularities) and/or model reduction. The SVD-based implementation methods possess a higher computational cost than the Cholesky-based algorithms because one SVD approximately costs two QR factorizations.

We start derivation of the SVD-based filtering methods with the prediction steps of the novel EM-0.5 DF-EKF and IT-1.5 DF-EKF methods in Algorithms~1 and~2. Again, we need to express formulas~\eqref{EM:Predict:P} and~\eqref{IT:Predict:P}, which are proved in Lemmas~1 and~2, in terms of the SVD factors $W_{P_{k|k-1}}$ and $\Sigma_{P_{k|k-1}}$, only. It is easy to show that for the EM-0.5 DF-EKF method we get
\begin{equation}
\biggl[W_{P_{k-1|k-1}^{(l+1)}},\Sigma^{1/2}_{P_{k-1|k-1}^{(l+1)}},(*) \biggr] \leftarrow {\mathbf svd}\biggl[\underbrace{\overline{{\mathbb F}{\mathbb X}}_{k-1|k-1}^{EM,\;(l)} \;\; \sqrt{\delta} G W_{Q}\Sigma_{Q}^{1/2}}_{pre-array \; A} \biggr] \label{eq:EM:SVD:P}
\end{equation}
where $W_{Q}$ and $\Sigma_{Q}^{1/2}$ are the SVD factors of the process covariance matrix $Q$, i.e. $Q = W_{Q}\Sigma_{Q} W_{Q}^{\top}$. The term $(*)$ denotes the orthogonal SVD factor of the pre-array, which is of no interest.

\begin{table*}[ht!]
{\small
\renewcommand{\arraystretch}{1.3}
\caption{The {\it SVD-based} continuous-discrete derivative-free EKF methods within the Euler-Maruyama discretization and It\^{o}-Taylor expansion.} \label{Tab:3}
\centering
\begin{tabular}{l||l|l}
\hline
& \cellcolor{myGray} {\bf Algorithm~1c: SVD-based EM-0.5 DF-EKF} & \cellcolor{myGray} {\bf Algorithm~2c: SVD-based IT-1.5 DF-EKF} \\
\hline
\hline
\textsc{Initialization:}  &  \multicolumn{2}{l}{0. SVD factorization $\Pi_0=W_{\Pi_0}\Sigma_{\Pi_0}W_{\Pi_0}^{\top}$. Set $\hat x_{0|0} = \bar x_0$, $W_{P_{0|0}} = W_{\Pi_0}$, $\Sigma_{P_{0|0}}^{1/2}=\Sigma_{\Pi_0}^{1/2}$ and $\alpha$; e.g. $\alpha = 1000$ suggested in~\cite{quine2006derivative}.}\\
\hline
 \textsc{Time} & \multicolumn{2}{l}{1. On $[t_{k-1}, t_k]$ introduce a mesh: $t_{k-1}^{(l)} = t_{k-1} + l \delta$, $l=0, \ldots, L-1$,  $\delta=\Delta/L$, $\Delta = |t_k-t_{k-1}|$, i.e. $t_{k-1}=t_{k-1}^{(0)} < \ldots t_{k-1}^{(l)}  < \ldots t_{k-1}^{(L)} = t_k$.} \\
\textsc{Update (TU):} & \multicolumn{2}{l}{\; Set up the initial values for integrators: $\hat x_{k-1|k-1}^{(0)}:=\hat x_{k-1|k-1}$ and $W_{P_{k-1|k-1}^{(0)}}:=W_{P_{k-1|k-1}}$, $\Sigma_{P_{k-1|k-1}^{(0)}}^{1/2}:=\Sigma_{P_{k-1|k-1}}^{1/2}$ at time node $t_{k-1}$.} \\
& \multicolumn{2}{l}{\; {\bf For} $l=1,\ldots L-1$ perform the following steps:}\\
& \multicolumn{2}{l}{\quad 2. Generate all sample points through ${\mathbb  X}_{k-1|k-1}^{(l)} = \hat x_{k-1|k-1}^{(l)} {\mathbf 1}^{\top} + \frac{\sqrt{n}}{\alpha} W_{P_{k-1|k-1}^{(l)}}\Sigma^{1/2}_{P_{k-1|k-1}^{(l)}}$.} \\
& \quad 3. Calculate estimate $\hat x_{k-1|k-1}^{(l+1)}$ by eq.~\eqref{EM:Predict:X}. & \quad 3. Calculate estimate $\hat x_{k-1|k-1}^{(l+1)}$ by eq.~\eqref{IT:Predict:X}. \\
 & \quad 4. Propagate sample points: ${\mathbb F}{\mathbb X}_{k-1|k-1}^{EM,\;(l)} =  f_d^{EM}\bigl(t_{k-1}^{(l)}, {\mathbb  X}_{k-1|k-1}^{(l)}\bigr)$. &
 \quad 4. Propagate ${\mathbb F}{\mathbb X}_{k-1|k-1}^{IT,\;(l)} =  f_d^{IT}\bigl(t_{k-1}^{(l)}, {\mathbb  X}_{k-1|k-1}^{(l)}\bigr)$.  \\
& \quad 5. Find $\overline{{\mathbb F}{\mathbb X}}_{k-1|k-1}^{EM,\;(l)} = \frac{\alpha}{\sqrt{n}} \left[{\mathbb F}{\mathbb X}_{k-1|k-1}^{EM,\;(l)} - f_d^{EM}\bigl(t_{k-1}^{(l)},\hat x_{k-1|k-1}^{(l)}\bigr)\;{\mathbf 1}^{\top} \right]$. & \quad 5. $\overline{{\mathbb F}{\mathbb X}}_{k-1|k-1}^{IT,\;(l)} = \frac{\alpha}{\sqrt{n}} \left[{\mathbb F}{\mathbb X}_{k-1|k-1}^{IT,\;(l)} - f_d^{IT}\bigl(t_{k-1}^{(l)},\hat x_{k-1|k-1}^{(l)}\bigr)\;{\mathbf 1}^{\top} \right]$.  \\
& \quad 6. Calculate SVD factors $W_{P_{k-1|k-1}^{(l+1)}}$, $\Sigma^{1/2}_{P_{k-1|k-1}^{(l+1)}}$ by eq.~\eqref{eq:EM:SVD:P}. &
\quad 6. Calculate SVD factors $W_{P_{k-1|k-1}^{(l+1)}}$, $\Sigma^{1/2}_{P_{k-1|k-1}^{(l+1)}}$ by eq.~\eqref{eq:IT:SVD:P}. \\
& \multicolumn{2}{l}{\; {\bf End}. At the last node, set $\hat x_{k|k-1}:=\hat x_{k-1|k-1}^{(L)}$ and $W_{P_{k|k-1}}:=W_{P_{k-1|k-1}^{(L)}}$, $\Sigma_{P_{k|k-1}}^{1/2}:=\Sigma_{P_{k-1|k-1}^{(L)}}^{1/2}$ at time node $t_{k-1}^{(L)}:=t_k$.}\\
\hline
\textsc{Measurement} & \multicolumn{2}{l}{7. Define sample points ${\mathbb  X}_{k|k-1}=\hat x_{k|k-1} {\mathbf 1}^{\top} + \frac{\sqrt{n}}{\alpha} W_{P_{k|k-1}}\Sigma_{P_{k|k-1}}^{1/2}$. Propagate $\hat z_{k|k-1} = h\left(k,\hat x_{k|k-1}\right)$ and ${\mathbb  Z}_{k|k-1} = h\left(k,{\mathbb  X}_{k|k-1}\right)$.} \\
\textsc{Update (MU):} & \multicolumn{2}{l}{8.  Define $P_{xz,k}=\overline{\mathbb X}_{k|k-1}\overline{\mathbb Z}_{k|k-1}^{\top}$ where $\overline{\mathbb  X}_{k|k-1} = \frac{\alpha}{\sqrt{n}} \bigl[{\mathbb  X}_{k|k-1} -  \hat x_{k|k-1}{\mathbf 1}^{\top} \bigr]$ and $\overline{\mathbb  Z}_{k|k-1}  =\frac{\alpha}{\sqrt{n}} \bigl[ {\mathbb  Z}_{k|k-1} -  \hat z_{k|k-1}{\mathbf 1}^{\top} \bigr]$.  } \\
& \multicolumn{2}{l}{9. Find $W_{R_{e,k}}$, $\Sigma_{R_{e,k}}^{1/2}$ by eq.~\eqref{SVD:Rek:1},  ${K}_{k}=P_{xz,k}W_{R_{e,k}}\Sigma_{R_{e,k}}^{-1}W_{R_{e,k}}^{\top}$, $\hat x_{k|k}=\hat x_{k|k-1}+{K}_k(z_k-\hat z_{k|k-1})$ and $W_{P_{k|k}}$, $\Sigma_{P_{k|k}}^{1/2}$ by eq.~\eqref{SVD:P:1}.} \\
\hline
\end{tabular}
}
\end{table*}

The equality above can be proved by multiplying the pre-array by its transpose and taking into account the properties of any orthogonal matrix, i.e. we have
\begin{align*}
AA^{\top} & = W_{P_{k-1|k-1}^{(l+1)}}\Sigma_{P_{k-1|k-1}^{(l+1)}}W_{P_{k-1|k-1}^{(l+1)}}^{\top} = P_{k-1|k-1}^{(l+1)} \\
  &  = \left[\overline{{\mathbb F}{\mathbb X}}_{k-1|k-1}^{EM,\;(l)} \quad \sqrt{\delta} G W_{Q}\Sigma_{Q}^{1/2}\right]
  \left[ \overline{{\mathbb F}{\mathbb X}}_{k-1|k-1}^{EM,\;(l)} \quad \sqrt{\delta} G W_{Q}\Sigma_{Q}^{1/2} \right]^{\top} \\
& =  \overline{{\mathbb F}{\mathbb X}}_{k-1|k-1}^{EM,\;(l)}\left(\overline{{\mathbb F}{\mathbb X}}_{k-1|k-1}^{EM,\;(l)}\right)^{\top} + \delta GQG^{\top},
\end{align*}
that is, formula~\eqref{EM:Predict:P} of the  EM-0.5 DF-EKF method.

Similarly, we re-derive formula~\eqref{IT:Predict:P} of the IT-1.5 DF-EKF in terms of the square-root factors of the predicted filter covariance matrix. For that, we take into account the re-arrangement of equation~\eqref{IT:Predict:P} in formula~\eqref{eq:proof:1} and conclude
{\small
\begin{align}
& \biggl[W_{P_{k-1|k-1}^{(l+1)}},\Sigma^{1/2}_{P_{k-1|k-1}^{(l+1)}},(*) \biggr]  \nonumber \\
&  \leftarrow {\mathbf svd} \underbrace{ \left[\overline{{\mathbb F}{\mathbb X}}_{k-1|k-1}^{IT,\;(l)}, \; \sqrt{\delta}(GW_{Q}\Sigma_{Q}^{1/2} + \frac{\delta}{2}{\mathbb L}f_{k-1}^{(l)}), \; \sqrt{\frac{\delta^3}{12}} {\mathbb L}f_{k-1}^{(l)} \right] }_{pre-array \; A} \label{eq:IT:SVD:P}
\end{align}
}
where ${\mathbb L}f_{k-1}^{(l)} := {\mathbb L} f\bigl(t_{k-1}^{(l)},\hat x_{k-1|k-1}^{(l)}\bigr)$. The equality above can be proved by multiplying the pre-array by its transpose, taking into account the properties of any orthogonal matrix and factorization in~\eqref{eq:proof:1}, and then comparing both sides with formula~\eqref{IT:Predict:P}.

Next, at the measurement update step, we need to compute the SVD factors of the residual covariance matrix that can be done in the following way:
\begin{equation}
\biggl[W_{R_{e,k}},\Sigma^{1/2}_{R_{e,k}},(*) \biggr] \leftarrow {\mathbf svd} \underbrace{\left[\overline{\mathbb Z}_{k|k-1} \quad W_{R_k}\Sigma_{R_k}^{1/2}\right]}_{pre-array \; A}
\label{SVD:Rek:1}
\end{equation}
where $W_{R_k}$ and $\Sigma_{R_k}^{1/2}$ are the SVD factors of the measurement covariance matrix $R_k$, i.e. $R_k = W_{R_k}\Sigma_{R_k} W_{R_k}^{\top}$. The term $(*)$ denotes the orthogonal SVD factor of the pre-array, which is of no interest.

Formula~\eqref{SVD:Rek:1} can be proved by multiplying the pre-array by its transpose and taking into account the properties of any orthogonal matrix, i.e. we have
\begin{align*}
AA^{\top} & = W_{R_{e,k}}\Sigma_{R_{e,k}}W_{R_{e,k}}^{\top} = R_{e,k} \\
& =  \left[\overline{\mathbb Z}_{k|k-1} \;\; W_{R_k}\Sigma_{R_k}^{1/2}\right]\left[\overline{\mathbb Z}_{k|k-1} \;\; W_{R_k}\Sigma_{R_k}^{1/2}\right]^{\top} = \overline{\mathbb Z}_{k|k-1}\overline{\mathbb Z}_{k|k-1}^{\top} + R_k,
\end{align*}
that is, exactly the first equation in~\eqref{eq:approx:Pxz:new} of the DF-EKF.

The formula for the filter's gain computation in terms of the SVD factors $W_{R_{e,k}}$ and $\Sigma_{R_{e,k}}$ is the following one:
\[
K_{k} = P_{xz,k}R_{e,k}^{-1} = P_{xz,k}W_{R_{e,k}}\Sigma_{R_{e,k}}^{-1}W_{R_{e,k}}^{\top}.
 \]

Finally, to calculate the SVD factors of the filtered covariance matrix, we define
{\small
\begin{equation}
\biggl[W_{P_{k|k}},\Sigma^{1/2}_{P_{k|k}},(*) \biggr] \leftarrow {\mathbf svd} \underbrace{\left[ (\overline{\mathbb X}_{k|k-1}-{K}_{k}\overline{\mathbb Z}_{k|k-1}) \quad K_k W_{R_k}\Sigma_{R_k}^{1/2}\right] }_{pre-array \; A}.  \label{SVD:P:1}
 \end{equation}
 }

Formula~\eqref{SVD:P:1} can be proved as follows:
\begin{align*}
AA^{\top} & = W_{P_{k|k}}\Sigma_{P_{k|k}}W_{P_{k|k}}^{\top} = P_{k|k} \\
 &= \left[\overline{\mathbb X}_{k|k-1}-{K}_{k}\overline{\mathbb Z}_{k|k-1}\right]\left[\overline{\mathbb X}_{k|k-1}-{K}_{k}\overline{\mathbb Z}_{k|k-1}\right]^{\top} + K_k R_k K_k^{\top}.
\end{align*}

To conclude this section, we summarize the SVD-based square-root methods derived for the EM-0.5 DF-EKF and IT-1.5 DF-EKF estimators in Table~\ref{Tab:3} in the form of Algorithms~1c and Algorithms~2c, respectively.

\section{Numerical experiments} \label{numerical:experiments}

To illustrate the performance of the novel derivative-free EKF methods, we provide a set of numerical tests on the radar tracking scenario from~\cite[Sec.~VIII]{Haykin2010} but with artificial ill-conditioned measurement scheme for provoking the numerical instability due to roundoff as discussed in~\cite[Example 7.1]{DyMc69} and~\cite[Section~7.2.2]{2015:Grewal:book}. Following~\cite{quine2006derivative}, all novel derivative-free EKF methods are implemented with $\alpha = 10^{3}$ that is shown to be sufficient for a convergence of the derivative-free EKF technique to the standard EKF methodology derived in~\cite{Ja70}. For a fair comparative study, the standard continuous-discrete EKF approach is implemented with the use of the Euler-Maruyama and It\^{o}-Taylor discretization schemes as well. One may easily obtain the related EKF implementation formulas, but for readers' convenience we refer to the summaries of the Euler-Maruyama-based EKF implementation method (EM-0.5 EKF) in~\cite[Algorithm~A1]{KuKu16SISCI} and the It\^{o}-Taylor-based EKF (IT-1.5 EKF) in~\cite[Algorithm~1]{KuKu18cIEEE_ICSTCC}. Finally, we additionally utilize the continuous-discrete CKF method proposed in~\cite{Haykin2010}, which has been designed within the It\^{o}-Taylor discretization scheme (i.e., the IT-1.5 CKF). The continuous-discrete CKF method based on the Euler-Maruyama discretization (i.e., the EM-0.5 CKF) can be also easily derived where the summary of computations is available in~\cite{KuKu20Automatica}.

\begin{exmp} \label{example:1:ill}
When performing a coordinated turn, the aircraft's dynamics obeys~\eqref{eq1.1} with the drift function
$f(\cdot)=\left[\dot{\epsilon}, -\omega \dot{\eta}, \dot{\eta}, \omega \dot{\epsilon}, \dot{\zeta},  0, 0\right]$ and the {\it standard} Brownian motion, i.e. $Q=I$ with $G={\rm diag}\left[0,\sigma_1,0,\sigma_1,0,\sigma_1,\sigma_2\right]$, $\sigma_1=\sqrt{0.2}$, $\sigma_2=0.007$.
The state consists of three positions, corresponding velocities and the turn rate, that is $x(t)= [\epsilon, \dot{\epsilon}, \eta, \dot{\eta}, \zeta, \dot{\zeta}, \omega]^{\top}$.  The turn rate is set to $\omega=3^\circ/\mbox{\rm s}$. The initial values are $\bar x_0=[1000\,\mbox{\rm m}, 0\,\mbox{\rm m/s}, 2650\,\mbox{\rm m},150\,\mbox{\rm m/s}, 200\,\mbox{\rm m}, 0\,\mbox{\rm m/s},\omega^\circ/\mbox{\rm s}]^{\top}$ and $P_0=I_7$. The state is observed through the measurement scheme
\begin{align*}
z_k & =
\begin{bmatrix}
1 & 1 & 1 & 1 & 1 &  1 &  1\\
1 & 1 & 1 & 1 & 1 &  1 &  1 + \gamma
\end{bmatrix}
x_k +
\begin{bmatrix}
v_k^1 \\
v_k^2
\end{bmatrix}, \; R = \gamma^{2}I_2
\end{align*}
where the ill-conditioning parameter $\gamma = 10^{-1},10^{-2},\ldots,10^{-14}$.
\end{exmp}

\begin{table*}
\caption{Degradation of accuracies ($\mbox{\rm ARMSE}$) of various Euler-Maruyama-based EKFs and CKF (with $L=512$ subintervals) while increasing ill-conditioning in Example~1.} \label{Tab:accuracy1}
\centering
{\small
\begin{tabular}{c||c|c||c|c|c|c|c}
\hline
& {\bf Standard EKF} & {\bf Cubature KF} & \multicolumn{5}{c}{\bf \texttt{EM-0.5 DF-EKF}: derivative-free EKF methods} \\
\cline{4-8}
$\gamma$ & \textit{Conventional}  & \textit{Conventional}  & \multicolumn{2}{c|}{Conventional implementations in \texttt{Alg.1}} & \multicolumn{3}{c}{Square-root implementations}  \\
\cline{4-8}
& (\texttt{EM-0.5 EKF})  & (\texttt{EM-0.5 CKF})  & {\scriptsize \it (with Cholesky in line 2,7)} &  {\scriptsize \it (with SVD in line 2,7)} & {\it Cholesky} \texttt{Alg.1a} & {\it Cholesky} \texttt{Alg.1b} & {\it SVD} \texttt{Alg.1c} \\
\hline
$10^{-1}$ & $4.376\cdot 10^{2}$ & $4.367\cdot 10^{2}$ &{\bf fail} & $4.352\cdot 10^{2}$ & $4.376\cdot 10^{2}$ & $4.376\cdot 10^{2}$ & $4.376\cdot 10^{2}$ \\
$10^{-2}$ & $4.375\cdot 10^{2}$ & {\bf fail} &&$4.439\cdot 10^{2}$ & $4.375\cdot 10^{2}$ & $4.375\cdot 10^{2}$ & $4.375\cdot 10^{2}$  \\
$10^{-3}$ & $4.375\cdot 10^{2}$ & &&$4.469\cdot 10^{2}$ & $4.375\cdot 10^{2}$ & $4.375\cdot 10^{2}$ & $4.375\cdot 10^{2}$ \\
$10^{-4}$ & $4.375\cdot 10^{2}$ & &&$4.332\cdot 10^{2}$ & $4.375\cdot 10^{2}$ & $4.375\cdot 10^{2}$ & $4.375\cdot 10^{2}$ \\
$10^{-5}$ & $4.375\cdot 10^{2}$ & && {\bf fail} & $4.375\cdot 10^{2}$ & $4.375\cdot 10^{2}$ & $4.375\cdot 10^{2}$  \\
$10^{-6}$ & $4.375\cdot 10^{2}$ & &&& $4.375\cdot 10^{2}$ & $4.375\cdot 10^{2}$ & $4.375\cdot 10^{2}$  \\
$10^{-7}$ & {\bf fail} & &&& $4.375\cdot 10^{2}$ & $4.375\cdot 10^{2}$ & $4.375\cdot 10^{2}$ \\
$10^{-8}$ &  & &&& $4.379\cdot 10^{2}$ & $4.374\cdot 10^{2}$ & $4.382\cdot 10^{2}$  \\
$10^{-9}$ &  & &&& $4.370\cdot 10^{2}$ & $4.236\cdot 10^{2}$ & $4.308\cdot 10^{2}$  \\
$10^{-10}$ & & &&& {\bf fail} & $3.721\cdot 10^{2}$ & {\bf fail}  \\
$10^{-11}$ &  & &&&  & $3.635\cdot 10^{2}$ &    \\
$10^{-12}$ &  & &&&  & $3.623\cdot 10^{2}$ &   \\
$10^{-13}$ &  & &&&  & $3.665\cdot 10^{2}$ &   \\
$10^{-14}$ &  & &&&  & $3.826\cdot 10^{2}$ &    \\
\hline
\end{tabular}
}
\end{table*}

To explore the numerical stability and the breakdown points of each filtering method under examination, the numerical experiments are provided for various values of $\gamma$. Given $P_0=I_7$, the residual covariance $R_{e,1} = HP_0H^{\top} + R$ and the observation matrix $H$ as well as $R$ are well-conditioned matrices when $\gamma$ is a large number. As $\gamma \to 0$, they become ill-conditioned matrices first and next $R_{e,1}$ (which needs to be inverted by the conventional filtering methods) becomes a singular one. In summary, our numerical experiment is intended on running several estimation scenarios starting from well-conditioned case and up to ill-conditioned situations in order to find the breakdown points of each filtering method under examination. In other words we observe the filter' divergence speed due to roundoff errors, which  characterizes their numerical instability.

In our numerical experiments, we solve a filtering problem on interval $[0s, 150s]$ with sampling period $\Delta = \Delta_k = 1s$ by fourteen EKF methods listed in Tables~\ref{Tab:accuracy1} and~\ref{Tab:accuracy2} for various ill-conditioned scenarios. More precisely, the set of numerical experiments are organized as follows. For each fixed value of the ill-conditioning parameter, $\gamma = 10^{-1},10^{-2},\ldots,10^{-14}$, we first solve the direct problem that is the numerical simulation of the given model. We discretize the stochastic system with a small stepsize $\delta = 0.0005$ on interval $[0s, 150s]$(s) to generate the {\it true} state vector $x^{\rm true}(t_k)$, $t_k \in [0s, 150s]$(s). Next, the measurement data is simulated with the sampling rates $\Delta = 1$(s). Given the data set, the inverse problem, i.e. the filtering problem, is solved to get the estimate of a hidden state, $\hat x_{k|k}$, over the time interval interval $[0, 150]$(s). We repeat the numerical test for $100$ trials and compute the accumulated root mean square error (ARMSE) by averaging over $100$ Monte Carlo runs and all seven entries of the state vector as follows:
 \begin{align}
\mbox{\rm ARMSE} & =\Bigl[\frac{1}{M
K}\sum_{M=1}^{100}\sum_{k=1}^K\sum_{j=1}^{n}\bigl(x^{\rm true}_{k,j}-\hat
x_{k|k,j}\bigr)^2\Bigr]^{1/2} \label{eq:acc:1}
\end{align}
where the subindex $j$, $j =1,\ldots,n$, refers to the $j$th entry of the $n$-dimensional state vector.

In each Monte Carlo run, we obtain its own simulated ``true'' state trajectory and the measurement data history. It is important that all estimators utilize the same initial conditions, the same simulated ``true'' state trajectory and the same measurements. The ARMSE values of the EM-0.5 EKF implementation methods are summarized in Table~\ref{Tab:accuracy1}. Meanwhile, the accuracies of the IT-1.5 EKF algorithms can be seen in Table~\ref{Tab:accuracy2}.

Having compared the ARMSEs in Tables~\ref{Tab:accuracy1} and~\ref{Tab:accuracy2}, we make the first conclusion. All Euler-Maruyama-based EKF methods are implemented with $L=512$ subintervals in the discretization mesh applied at the prediction filtering step, meanwhile the It\^{o}-Taylor-based algorithms are implemented with much less subdivision steps, which is $L=64$. As can be seen, although the  Euler-Maruyama-based EKF methods are implemented with considerably smaller discretization step $\delta$ (because $L=512$ is larger than $L=64$) the It\^{o}-Taylor-based EKF algorithms outperform the Euler-Maruyama-based EKF methods for estimation accuracies, significantly. This result is inline with the theory of numerical methods since the It\^{o}-Taylor expansion is a higher order method than the Euler-Maruyama scheme.

We conclude that the estimation methods designed with the use of It\^{o}-Taylor expansion, i.e. in our case any IT-1.5 EKF-type implementation in Table~\ref{Tab:accuracy2}, are more accurate and they require less subdivisions in the underlying discretization scheme to provide sufficient estimation quality compared to their EM-0.5 EKF counterparts in Table~\ref{Tab:accuracy1}. However, the IT-1.5 EKFs are rather complicated filtering methods and, hence, the EM-0.5 EKF-type estimators still may have some merit for solving practical applications. Indeed, they are easy to implement and they provide an adequate accuracy when the number of subdivision steps is high enough to maintain a good estimation quality. These are the attractive features of the Euler-Maruyama-based filtering methods.

\begin{table*}
\caption{Degradation of accuracies ($\mbox{\rm ARMSE}$) of various It\^{o}-Taylor-based EKFs and CKF (with $L=64$ subintervals) while increasing ill-conditioning in Example~1.} \label{Tab:accuracy2}
\centering
{\small
\begin{tabular}{c||c|c||c|c|c|c|c}
\hline
& {\bf Standard EKF} & {\bf Cubature KF} &  \multicolumn{5}{c}{\bf \texttt{IT-1.5 DF-EKF}: derivative-free EKF methods} \\
\cline{4-8}
$\gamma$ & \textit{Conventional}  & \textit{Conventional}  & \multicolumn{2}{c|}{Conventional implementations in \texttt{Alg.2}} & \multicolumn{3}{c}{Square-root implementations}  \\
\cline{4-8}
&  (\texttt{IT-1.5 EKF})  & (\texttt{IT-1.5 CKF})  & {\scriptsize \it (with Cholesky in line 2,7)} &  {\scriptsize \it (with SVD in line 2,7)} & {\it Cholesky} \texttt{Alg.2a} & {\it Cholesky} \texttt{Alg.2b} & {\it SVD} \texttt{Alg.2c} \\
\hline
$10^{-1}$ & $3.743\cdot 10^{1}$ & $1.726\cdot 10^{1}$ & {\bf fail} & $4.407\cdot 10^{1}$ & $1.055\cdot 10^{1}$ & $1.055\cdot 10^{1}$ & $1.055\cdot 10^{1}$ \\
$10^{-2}$ & $4.433\cdot 10^{1}$ & {\bf fail}  &&$3.368\cdot 10^{1}$ & $1.058\cdot 10^{1}$ & $1.058\cdot 10^{1}$ & $1.058\cdot 10^{1}$ \\
$10^{-3}$ & $4.138\cdot 10^{1}$ &&&$3.806\cdot 10^{1}$ & $1.058\cdot 10^{1}$ & $1.058\cdot 10^{1}$ & $1.058\cdot 10^{1}$ \\
$10^{-4}$ & $4.139\cdot 10^{1}$ &&&$3.447\cdot 10^{1}$ & $1.058\cdot 10^{1}$ & $1.058\cdot 10^{1}$ & $1.058\cdot 10^{1}$ \\
$10^{-5}$ & $4.138\cdot 10^{1}$ &&&{\bf fail} & $1.058\cdot 10^{1}$ & $1.058\cdot 10^{1}$ & $1.058\cdot 10^{1}$ \\
$10^{-6}$ & $4.138\cdot 10^{1}$ &&& & $1.058\cdot 10^{1}$ & $1.058\cdot 10^{1}$ & $1.058\cdot 10^{1}$ \\
$10^{-7}$ & {\bf fail}  &&& & $1.055\cdot 10^{1}$ & $1.058\cdot 10^{1}$ & $1.058\cdot 10^{1}$ \\
$10^{-8}$ &  &&& & $1.055\cdot 10^{1}$ & $1.059\cdot 10^{1}$ & $1.046\cdot 10^{1}$ \\
$10^{-9}$ &  &&& & $1.911\cdot 10^{1}$ & $1.103\cdot 10^{1}$ & $2.551\cdot 10^{1}$ \\
$10^{-10}$ &  &&&  & {\bf fail} & $1.159\cdot 10^{1}$ & {\bf fail} \\
$10^{-11}$ &  &&&  &  & $3.106\cdot 10^{1}$ &  \\
$10^{-12}$ &  &&& &  & $3.475\cdot 10^{1}$ &  \\
$10^{-13}$ &  &&&  &  & $2.472\cdot 10^{1}$ &  \\
$10^{-14}$ &  &&&  &  & $3.015\cdot 10^{1}$ &  \\
\hline
\end{tabular}
}
\end{table*}

\begin{figure}
\includegraphics[width=0.5\textwidth]{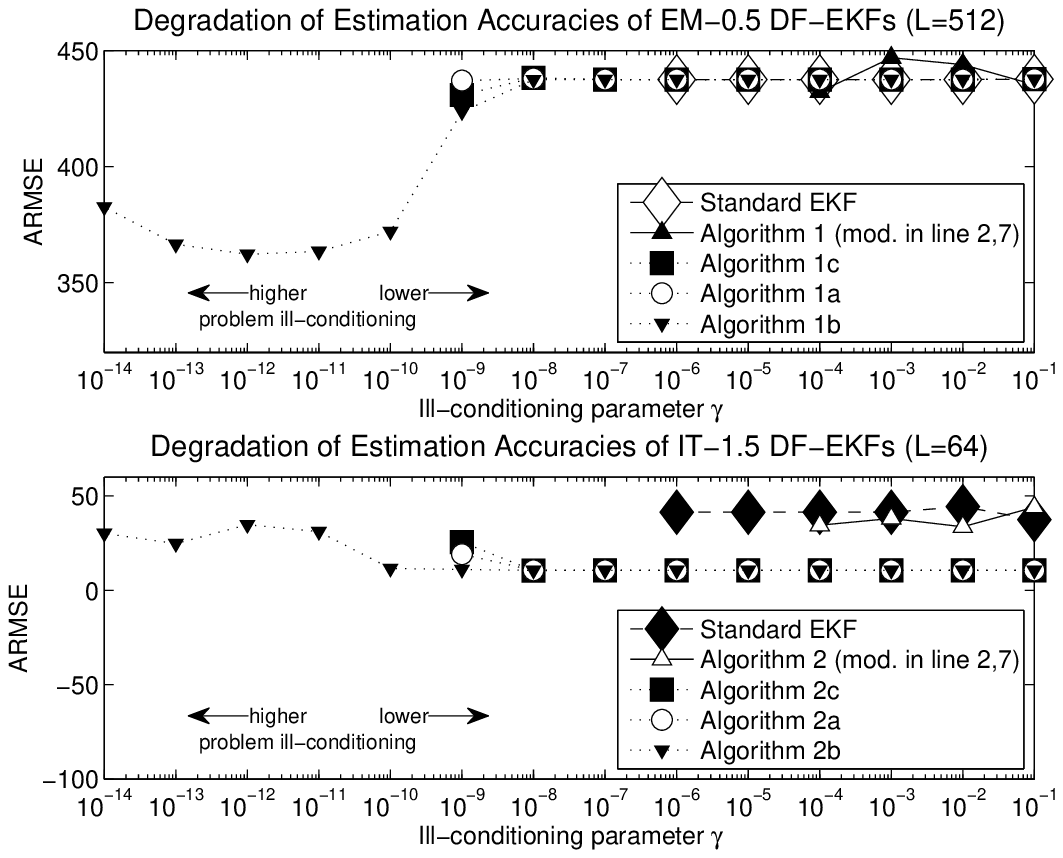}
\caption{The accuracies  degradation of various continuous-discrete EKF estimators on the ill-conditioned problem in Example~\ref{example:1:ill}.} \label{fig:1}
\end{figure}

The most significant finding of our numerical tests concerns the numerical stability of the filtering methods derived. Indeed, the key goal of Example~1 is to highlight some insights about the robustness of the filtering methods under examination with respect to roundoff errors when the test problems become ill-conditioned. More precisely, we should pay an attention to the {\it breakdown} points of each filtering method when the ill-conditioning parameter $\gamma$ tends to a machine precision limit. These points are marked by `{\bf fail}' in Tables~\ref{Tab:accuracy1} and~\ref{Tab:accuracy2} and this means that the particular numerical method fails in real computations. We additionally illustrate the results by Fig.~\ref{fig:1} where the EM-0.5 EKF and IT-1.5 EKF implementation methods are plotted separately for a better exposition.

To start our discussion, we compare the results obtained by the standard EKF approach and the derivative-free EKF methodology, first. For that, one should compare the results of the conventional implementation methods summarized in the second column (see \texttt{Standard EKF}) and the fourth column (see \texttt{EM-0.5 DF-EKF} in Algorithm~1) in Table~\ref{Tab:accuracy1}. It is clearly seen that the novel derivative-free EKF method is unstable with respect to roundoff errors compared to  its standard EKF counterpart because the novel Algorithm~1 fails for ill-conditioning parameter $\gamma = 10^{-1}$ but the breakdown point of the standard EKF is $\gamma = 10^{-7}$. The same conclusion holds true for the IT-1.5 EKF methods in Table~\ref{Tab:accuracy2}. Again, the conventional \texttt{IT-1.5 DF-EKF} in Algorithm~2 is the most vulnerable algorithm to roundoff errors in Example~1. Such numerical behavior of the conventional derivative-free EKF methods has been anticipated in Section~\ref{sec:main:conventional} when  novel  Algorithms~1 and~2 were discussed. The key problem of their numerical instability is the requirement of Cholesky decomposition at each iterate of the derivative-free filtering methods to generate the sample vectors; see lines~2 and~7 in both conventional Algorithms~1 and~2. Example~1 is designed in such a way that it highlights the insights of numerical instability problem. The roundoff errors in the scenarios of Example~1 destroy the theoretical properties of the filter covariance matrices that makes the Cholesky decomposition unfeasible. This interrupts the practical calculations and Algorithms~1, 2 fail. Recall, the standard EKF methodology does not require any matrix factorization and the sample points generation and, as a result, the standard EKF is more stable with respect to roundoff errors. This explains a better numerical robustness of the standard EKF framework compared to the derivative-free EKF estimation strategy.

Besides, the standard EKF methods are more stable with respect to roundoff errors than the CKF algorithms under examination. It is clearly seen from the results presented in the second columns (see \texttt{Standard EKF}) and the third columns (see \texttt{Cubature KF} algorithms) in Tables~\ref{Tab:accuracy1} and \ref{Tab:accuracy2}. Indeed, the EM-0.5 CKF and IT-1.5 CKF methods fail when $\gamma < 10^{-1}$, meanwhile, the breakdown point of the standard EKF is $\gamma = 10^{-7}$. Additionally, we note that the CKF algorithms are slightly more stable than the DF-EKF methods for this type of ill-conditioned tests. However, in general, the EM-0.5 CKF and IT-1.5 CKF methods inherit the same numerical instability problem as all derivative-free filtering methods when the matrix factorization is requested at each iterate of the filter in order to generate the sigma/cubature/quadrature vectors. Similarly to the DF-EKF methodology, the CKF algorithms fail because the roundoff errors  destroy the theoretical properties of the filter covariance matrices that makes the Cholesky decomposition unfeasible.

As discussed in Section~\ref{sec:main:conventional}, we may replace Cholesky decomposition in lines~2 and~7 in conventional Algorithms~1 and~2 by SVD factorization to generate the sample points. Having analyzed the obtained results of the conventional implementations summarized in the fifth columns in Tables~\ref{Tab:accuracy1} and~\ref{Tab:accuracy2}, we conclude that such refined modification of the conventional algorithms certainly improves their numerical stability with respect to roundoff errors. Indeed, the breakdown point of such {\it conventional} derivative-free EKF implementations is $\gamma = 10^{-5}$, i.e. they work accurately and sustainedly in the well-conditioned and mild ill-conditioned scenarios in Example~1.

Additionally, we remark that the estimation accuracies of the novel derivative-free EKF estimators and their standard EKF counterparts as well as the Cubature KFs in Tables~\ref{Tab:accuracy1} and~\ref{Tab:accuracy2}, are similar, but the CKF and DF-EKF methods are slightly more accurate. Our results are also in line with the conclusion made in~\cite{quine2006derivative} where $\alpha = 1000$ is shown to be enough to ensure an excellent convergence of the derivative-free EKF to the standard EKF. We conclude that the novel derivative-free EKF methods work with a good estimation quality on the well-conditioned problems and slightly outperform the standard EKF technique for estimation quality, but the derivative-free EKF framework is vulnerable to roundoff errors. Thus, the derivation and practical utilization of the square-root implementation methods in case of the derivative-free EKF framework is an extremely important issue.

To strengthen the numerical robustness of the novel derivative-free EKF methods, the square-root implementations derived in Section~\ref{sec:main:squareroot} should be applied. We next examine their numerical robustness by comparing the results summarized in the last three columns in both Tables~\ref{Tab:accuracy1} and~\ref{Tab:accuracy2}. Having compared the breakdown points of the square-root algorithms proposed, we conclude that the Cholesky-based derivative-free EKF Algorithms~1b and~2b are the most stable methods in a finite precision arithmetics. Indeed, their Cholesky-based counterparts, i.e. the \texttt{Cholesky EM-0.5 DF-EKF} in Algorithm~1a and the \texttt{Cholesky IT-1.5 DF-EKF} in Algorithm~2a, fail for ill-conditioned tests with $\gamma = 10^{-10}$ as well as their SVD-based counterparts, i.e. the \texttt{SVD EM-0.5 DF-EKF} in Algorithm~1c and the \texttt{SVD IT-1.5 DF-EKF} in Algorithm~2c. Meanwhile, the \texttt{Cholesky EM-0.5 DF-EKF} in Algorithm~1b and the \texttt{Cholesky IT-1.5 DF-EKF} in Algorithm~2b work in a stable way, i.e. without a failure, for all ill-conditioned tests under examination. Recall, two types of the Cholesky-based square-root DF-EKF implementations differ by the number of QR factorizations implemented for computing the square-root factors at the measurement update steps. This fact also has a strong impact on the gain matrix computation. More precisely, Algorithms~1a and~2a utilizes two QR factorizations at the measurement update steps and, next, the square-root factors $R_{e,k}^{1/2}$ and $R_{e,k}^{\top/2}$ should be inverted to compute the gain matrix $K_k$ in line~9 of Algorithms~1a and~2a, respectively. Meanwhile, Cholesky-based Algorithms~1b and~2b imply only one QR factorization, which is implemented to the unique pre-array at the measurement update steps. This allows for calculating the normalized cross-covariance $\bar P_{xz}$, which is simply read-off from the post-array after the orthogonal rotation. The accessibility of this term yields one less matrix inversion (which is $R_{e,k}^{-1/2}$, only) for calculating the gain matrix $K_k$ in line~9 of Algorithms~1b and~2b, respectively. Finally, the SVD factorization-based Algorithms~1c and~2c possess the same numerical stability as the Cholesky-based Algorithms~1a and~2a.

\begin{figure}
\includegraphics[width=0.5\textwidth]{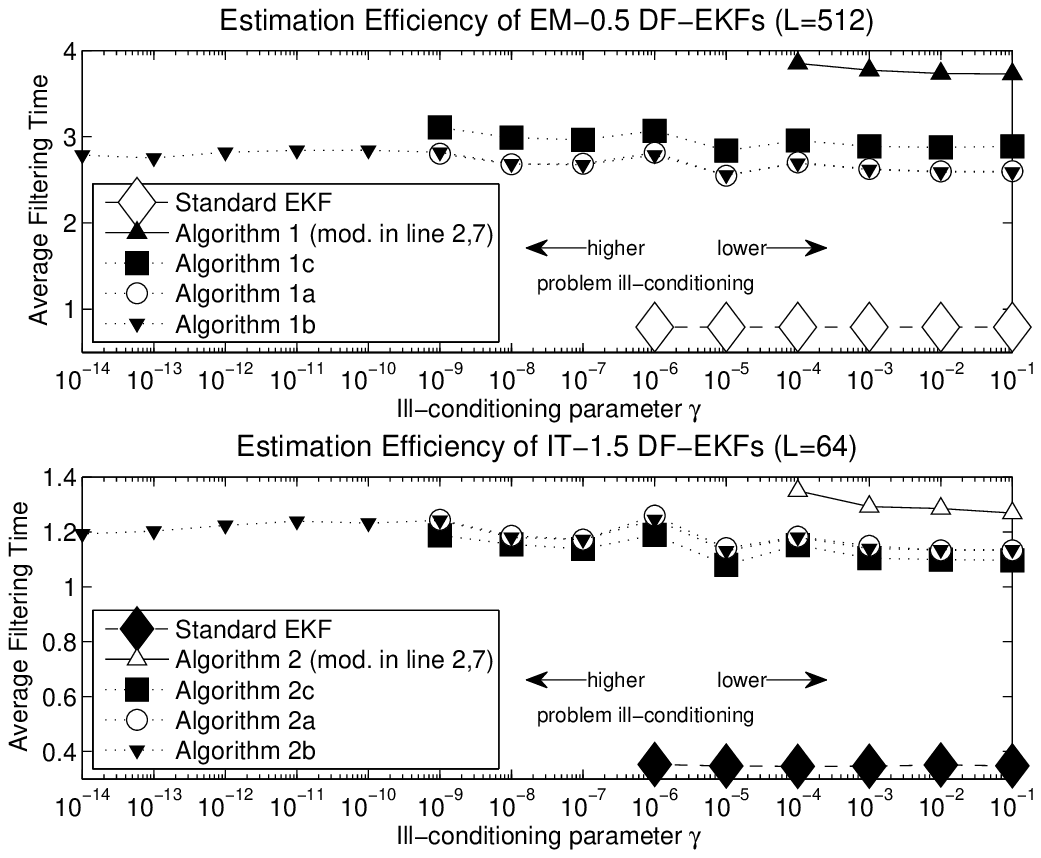}
\caption{Estimation Efficiency in Run Time (s) of various continuous-discrete EKF estimators on the ill-conditioned problem in Example~\ref{example:1:ill}.} \label{fig:2}
\end{figure}

In our last set of numerical tests, we examine an efficiency of the methods proposed. For that, we additionally compute the average CPU time (in seconds) over $100$ Monte Carlo runs for each estimator under discussion and illustrate the resulted values by Fig.~\ref{fig:2}. The results of the EM-0.5 EKF and IT-1.5 EKF implementation methods are plotted separately for a better exposition and fair comparative study. As can be seen, the standard EKF techniques are faster than the related derivative-free EKF variants. This result has been anticipated since the DF-EKF methodology requires the sample vectors generation, their propagation and matrix factorizations, in contrast to the standard EKF without such extra computations. Next, we may conclude that the CPU time of the Cholesky-based DF-EKF implementations are almost the same. To  observe this fact, one should compare the results of Algorithm~1a with the outcomes of Algorithm~1b as well as the results of Algorithm~2a with the outcomes of Algorithm~2b. Recall, the difference in these implementation methods is one extra QR factorization in Algorithms~1a and~2a, compared to  Algorithms~1b and~2b. Additionally, the SVD-based Algorithms~1c and~2c are slightly slower than other square-root DF-EKF implementation methods. Finally, the conventional DF-EKF implementations with the SVD utilized for the sample vectors generation in lines~2 and~7 of Algorithms~1 and~2 are the slowest implementations although the difference with other DF-EKF algorithms is not substantial.

\section{Conclusion}\label{Section:conclusion}

In this paper, we developed the {\it continuous-discrete} derivative-free EKF methods within the Euler-Maruyama and It\^{o}-Taylor discretization schemes applied for the stochastic nonlinear systems. The results of numerical tests substantiate a high estimation quality of the novel derivative-free EKF methods and, additionally, show that they either equal or slightly outperform the standard EKF techniques for accuracy. However, the derivative-free EKF framework is vulnerable to roundoff errors. Being a higher order method, the It\^{o}-Taylor expansion yields the related It\^{o}-Taylor-based DF-EKF algorithms that maintain a better estimation quality with a less number of subdivisions in the discretization mesh to be implemented than the Euler-Maruyama-based DF-EKF alternatives. Meanwhile, a simplicity and good filtering accuracy when the number of subdivision steps is high enough are still the attractive features of the Euler-Maruyama-based filtering methods.

We also stress that the novel filters are derivative-free methods and, hence, they are especially effective for working with stochastic systems with highly nonlinear and/or nondifferentiable drift and observation functions, i.e. when the calculation of Jacobian matrices are either problematical or questionable. The investigation of the novel filtering methods on the ill-conditioned test problems has shown that the price to be paid for these benefits is the degraded numerical stability of the conventional implementations with respect to roundoff errors. To resolve this problem, we have additionally derived the stable square-root methods within both the Cholesky and SVD square-root factorizations. One of the most significant findings to emerge from this study is that the novel square-root continuous-discrete EKF methods are numerically stable to roundoff and they possess all benefits of the derivative-free estimators. Thus, the practical utilization of the square-root implementation methods in case of the derivative-free EKF framework is an extremely important issue. We may conclude that the results of this study will have a number of important implications for future EKF application for solving practical problems.

\section*{Acknowledgements}
The authors acknowledge the financial support of the Portuguese FCT~--- \emph{Funda\c{c}\~ao para a Ci\^encia e a Tecnologia},
through the \emph{Scientific Employment Stimulus - 4th Edition} (CEEC-IND-4th edition) programme (grant number 2021.01450.CEECIND) and through the projects UIDB/04621/2020 and UIDP/04621/2020 of CEMAT/IST-ID, Center for Computational and Stochastic Mathematics, Instituto Superior T\'ecnico, University of Lisbon.

\section*{References}
\bibliographystyle{model1b-num-names}

\end{document}